\documentclass[a4paper,11pt]{amsart}

\usepackage[hidelinks]{hyperref}

\usepackage{amsmath,amsthm,amssymb,latexsym,amsfonts,xcolor}
\usepackage[utf8]{inputenc} % for older LaTeX; newer LaTeX defaults to UTF-8
\usepackage[english,activeacute]{babel}
\usepackage{graphicx}

\usepackage{enumerate}
\usepackage{csquotes}

\usepackage[a4paper, margin=1.8cm]{geometry}
\usepackage{parskip}
\usepackage{hyperref}
\usepackage{cleveref}

\usepackage{cancel}
\usepackage[normalem]{ulem}

\allowdisplaybreaks

\theoremstyle{plain}

\newtheorem*{thm*}{Theorem}
\newtheorem{thm}{Theorem}[section]
\newtheorem{conj}[thm]{Conjecture}

\newtheorem{rk}[thm]{Remark}

\theoremstyle{definition}

\theoremstyle{remark}

\theoremstyle{plain}

\newcommand{\sinc}{\mathrm{sinc}}
\newcommand{\R}{\mathbb{R}}
\newcommand{\T}{\mathbb{T}}
\newcommand{\D}[1]{D_x^{#1}}
\newcommand{\Hs}[2][\beta]{\mathcal{H}_{s,#1}(\T ^{#2})   }

\newcommand{\SN}{   \mathcal{S}_{\alpha,\beta}^N(t)}
\newcommand{\NN}{   \mathcal{N}_{\alpha,\beta}^N(t)}

\newcommand{\Sinf}{   \mathcal{S}_{\alpha,\beta}^{N,\infty}}
\newcommand{\Ninf}{   \mathcal{N}_{\alpha,\beta}^{N,\infty}}
\title[Numerical study of probabilistic well-posedness]{Numerical study of probabilistic well-posedness of one dimensional fractional nonlinear wave equations}

\author[W. Ruffenach]{Wandrille Ruffenach}
\address{ENS de Lyon, CNRS, LPENSL, UMR5672, 69342, Lyon cedex 07, France.}

\author[N. Tzvetkov]{Nikolay Tzvetkov}
\address{ENS de Lyon, CNRS, UMPA, 69342, Lyon cedex 07, France.}

\begin{document}
	
\begin{abstract}
The three dimensional cubic defocusing nonlinear wave equation is known to be ill-posed for general low regularity initial data. Well-posedness can however be recovered globally in time on a probabilistic level, provided the Fourier coefficients of the random initial data follow a sub-Gaussian law, an assumption on which the available proofs rely heavily. In this article we perform numerical simulations of the one dimensional fractional cubic defocusing wave equation in a periodic setting, with low regularity random initial data drawn either from a Gaussian or from a heavy tailed Student law. Besides illustrating probabilistic well-posedness and norm inflation in both the energy subcritical and supercritical regimes, our simulations display no qualitative difference between the heavy tailed and the sub-Gaussian cases. This leads us to conjecture that probabilistic well-posedness holds under the sole assumption that the Fourier coefficients are centered with finite variance.
\end{abstract}

	\maketitle

\section{Introduction}\label{sec:Intro}

We consider the Cauchy problem for the defocusing fractional nonlinear wave equation on the torus $\T^d= \left( \R / \mathbb{Z}\right)^d $
\begin{equation}
	\begin{cases}
		&\partial^2_t u +\D{2\beta} u+u^3=0, \, x \in \T^d \\
		&(u(0,x), \partial_t u(0,x))=(u_0(x),\dot{u}_0(x) )\in H^s(\T^d)\times H^{s-\beta}(\T^d),
	\end{cases}
	\label{eq:NLWaves} \tag{fNLW}
\end{equation}
where $\D{}=\sqrt{-\Delta}$ is the Fourier multiplier by $|2\pi k|$, $\beta>0$ and for $s\in \R$, the Sobolev space $H^s(\T^d)$ is
the completion of $C^\infty(\T^d)$ with respect to the norm
$$ 
\|v\|^2_{ H^s(\T^d)}=\sum_{n \in \mathbb{Z}^d} \left\langle n \right\rangle^{2s} |\widehat{v}_n |^2\, ,
%H^s(\T^d) = \Big \lbrace \left. v \in L^2(\T^d)\,\, \right| \,\, \sum_{n \in \mathbb{Z}^d} \langle n \rangle^{2s} |\widehat{v}_n |^2 <+\infty   \Big \rbrace, 
$$
where $ \left\langle n \right\rangle= \sqrt{(2\pi n )^2 +1}$ and $\widehat{v}_n=   \int_{[0,1]^d} v(x) e^{-2i\pi  n \cdot x} dx   $ are the Fourier modes of $v$. In addition, we introduce the space
$$ \Hs{d}  = H^s(\T^d) \times H^{s-\beta}(\T^d)  $$
endowed with the norm $$ \|(u,v)\|_{\Hs{d}} =  \|u\|_{H^{s}(\T^{d})} + \|v\|_{H^{s-\beta}(\T^{d})}.$$
The equation \eqref{eq:NLWaves} enjoys a Hamiltonian structure with the  conservation of the Hamiltonian
\begin{equation} \label{eq:HamiltonianNLWave}
	\mathcal{H}[u]= \int \Big[\left( \partial_t u \right)^2+(\D{\beta} u)^2+\dfrac{1}{2}u^4\Big] \, \mathrm{d}x,
\end{equation}
by the flow of \eqref{eq:NLWaves} provided $u$ is regular enough. The space $\Hs{}$ appears to be a natural framework for the well-posedness theory of \eqref{eq:NLWaves} in view of the structure of the Hamiltonian \eqref{eq:HamiltonianNLWave}. In addition, if $u$ is a solution of \eqref{eq:NLWaves} and $\lambda >0$ so does
\begin{equation}
	u_\lambda(t,x)=\lambda^\beta u(\lambda^\beta t , \lambda x) .
\end{equation}
In particular, if $u$ solves \eqref{eq:NLWaves} with a lifespan $T$, then $u_\lambda$ solves \eqref{eq:NLWaves} on a rescaled torus with a lifespan $T_\lambda=\lambda^{-\beta}T$ and 
\begin{equation}
	\begin{cases}
		&  \|(u_\lambda(0,\cdot),\partial_t u_\lambda(0,\cdot)) \|_{\Hs{d}  } \simeq \lambda^{s-(\frac d2-\beta) }  \|(u_0, \dot{u}_0) \|_{\Hs{d} },\\
		&  \mathcal{H}[u_\lambda(0, \cdot)]= \lambda^{4\beta -d }  \mathcal{H}[u(0, \cdot)]\, .
	\end{cases}
\end{equation}
This scaling symmetry introduces a critical regularity $s_c= \dfrac{d}{2}-\beta$ and a critical dispersion parameter $\beta_c= d/4$. In the following, we will refer to the case $\beta> \beta_c$ as energy subcritical and $\beta< \beta_c$ as energy super-critical.

Note that nonlinear fractional wave equations such as \eqref{eq:NLWaves} appears in several physically relevant systems such as the long-range $\phi^4$ model at criticality \cite{tara2006fractional}.

In the present article, we will study numerically the local well-posedness of \eqref{eq:NLWaves} in the one dimensional case $d=1$ for various dispersion parameters $\beta$. The restriction to one space dimension allows to perform highly resolved numerical simulations for a reasonable cost. Then, changing the dispersion relation set by $\beta$ allows us to explore both the energy subcritical  and supercritical regimes.
In the regime investigated in this paper the Knapp obstruction to deterministic well-posedness does not seem present. It is known that such an obstruction is present in the two dimensional case for $\beta=1$. We refer to the book by Tao \cite{tao2006nonlinear} for details on the Knapp obstruction.
%%%%%%%%%%%%%%%%%%%%%%%%
\subsection{Deterministic well-posedness and ill-posedness}\label{sec:LocalWP}
Existing results regarding the well-posedness of the cubic defocusing wave equation are mainly established in the energy-subcritical case for which, roughly speaking, the dispersion is stronger than the nonlinearity. Let us now discuss local and global well-posedness in the energy subcritical case with regular enough ($s>s_c$) initial data and present the existing results of ill-posedness and probabilistic well-posedness in the low regularity ($s<s_c$) energy subcritical regime. 
\subsubsection{On well-posedness in the energy subcritical regime}
Following \cite{Tzvetkov2019random}, we will sketch the proof of local and global in time well-posedness in the energy subcritical case for $\beta \geq d/3 > \beta_c=d/4$ and $s \geq \beta > s_c$ (i.e. $s$ is allowed to go only down to $\beta$ which in the energy subcritical case $\beta>d/4$ is strictly larger than $s_c= d/2-\beta$). For $S_{u_0,\dot{u}_0}$ the flow of the linear part of \eqref{eq:NLWaves} and $\Phi_{u_0,\dot{u}_0}$ the flow of \eqref{eq:NLWaves} Duhamel's formulation of \eqref{eq:NLWaves} writes as follows 
\begin{equation}
	\Phi_{u_0,\dot{u}_0}(t)= S_{u_0,\dot{u}_0}(t) - \int_0^t \dfrac{\sin \left((t-\tau) \D{\beta} \right)}{\D{\beta}} (\Phi_{u_0,\dot{u}_0}(\tau))^3 \, d\tau.
\end{equation}
Then, with the Banach space $X_T= C\Big([0,T],H^{\beta}(\T^d ) \Big)$ and the norm $ \|u\|_{X_T}= \sup_{t \in[0,T]} \|u\|_{H^\beta(\T^d)}$ it is easy to check that 
$$ \|S_{u_0,\dot{u}_0}\|_{X_T} \lesssim \| u_0 \|_{H^\beta(\T^d)} + \| \dot{u}_0 \|_{L^2(\T^d)}.  $$ 
Using this estimate and the Sobolev embedding $H^{\beta}(\T^d) \subset L^6(\T^d)$ for $\beta \geq d/3$ we obtain the following bound for the nonlinear flow,
\begin{align*}
	\|\Phi_{u_0,\dot{u}_0} \|_{X_T} \leq& C \Big( \| u_0 \|_{H^\beta(\T^d)} + \| \dot{u}_0 \|_{L^2(\T^d)}+ T \sup_{t\in[0,T]} \|\D{-\beta }\left( \Phi_{u_0,\dot{u}_0}\right)^3 \|_{H^{\beta}(\T^d)}  \Big) \\
	&\leq  C \Big( \| u_0 \|_{H^\beta(\T^d)} + \| \dot{u}_0 \|_{L^2(\T^d)}+ T \sup_{t\in[0,T]} \|\Phi_{u_0,\dot{u}_0}\|^3_{L^6(\T^d)}  \Big)  \\
	&\leq C \Big( \| u_0 \|_{H^\beta(\T^d)} + \| \dot{u}_0 \|_{L^2(\T^d)}+ T  \| \Phi_{u_0,\dot{u}_0}\|^3_{X_T}  \Big).  
\end{align*}
Using this bound in a fixed point argument in $X_T$ yield local well-posedness in $H^{\beta}(\T^d) \times L^2(\T^d)$ for $\beta \geq \frac d3$. Now, if the initial condition lies in $\Hs{}$ for some $s\geq \beta \geq \frac d3$, the regularity of the initial condition will propagate to later times. This propagation of regularity follows from Kato-Ponce inequality and the Sobolev embedding $H^{\beta}(\T^d) \subset L^6(\T^d)$ for $\beta \geq \frac d3$. In particular, they are applied to obtain the bound 
\begin{align*}
	\|\D{-\beta}\left( \Phi_{u_0,\dot{u}_0}(\tau)\right)^3 \|_{H^{s}(\T^d)}  &\lesssim  \|\D{s-\beta} \Phi_{u_0,\dot{u}_0}(\tau)\|_{L^6(\T^d)}\| \Phi_{u_0,\dot{u}_0}(\tau)\|^2_{L^6(\T^d)},\\&  \lesssim \|\D{s} \Phi_{u_0,\dot{u}_0}(\tau)\|_{L^2(\T^d)}   \| \Phi_{u_0,\dot{u}_0}(\tau)\|^2_{H^{\beta}(\T^d)} ,
\end{align*}
which ultimately yield propagation of regularity after following the same steps for $\partial_tu$. Then, using the conservation of energy \eqref{eq:HamiltonianNLWave}, it is possible to extend to global well-posedness of \eqref{eq:NLWaves} in $\Hs{d}$ for $s\geq \beta$ and $\beta \geq d/3$. 
%Using Strichartz's estimates, we believe that it is possible to establish the same result for $\beta > \beta_c= \frac d4$. 

For lower regularity $s\in (s_c,\beta)$, we expect local well-posedness to hold which should follow from Strichartz's estimates and we conjecture it can be extended to global well posedness. A result going in that direction can be established for the energy subcritical ($\beta=1$, $d=3$) case where global well-posedness holds for any $ \frac12 =s_c < \frac{13}{18} <s< \beta$, see \cite{Tzvetkov2019random} and references therein.  

\subsubsection{On ill-posedness}
For initial data with (positive) regularity lower than the critical regularity, local well-posedness does not hold in general. Based on ideas of Lebeau \cite{lebeau2005} one can prove the following result for $\beta=1$ and $d=3$.
\begin{thm}[\protect{\cite[Theorem 1.33]{Tzvetkov2019random}}]	\label{thm:IllPosedness} 
	Let $s\in(0,1/2)$ and $(u_0,\dot{u}_0) \in \Hs[1]{3}$. 
	There exists a sequence 
	$$ u_N \in C^0(\R; C^\infty(\T^3)), \quad N=1,2, \dots $$
	such that 
	$$(\partial^2_t -\Delta) u_N +u_N^3=0, $$
	with 
	$$ \lim_{N \rightarrow \infty } \| (u_N(0)-u_0, \partial_t u_N(0) -\dot{u}_0  )\|_{\Hs[1]{3}}  =0$$
	but for all $T>0$,
	$$ \lim_{N \rightarrow \infty } \| (u_N(t), \partial_t u_N(t))   \|_{L^\infty\left( [0,T], \Hs[1]{3}\right)}  =+\infty.$$
\end{thm}
This shows that for low regularity initial data, the non-linear wave equation lacks continuous dependence with respect to the initial data. In the numerical study performed in this paper, we will focus on the one dimensional problem on the torus in both energy subcritical and energy super-critical cases. Even though \eqref{eq:NLWaves}  with $\beta=1$ is ill-posed in general in $\Hs[1]{3}$ for any $s<1/2$, it is possible to restore well-posedness on a probabilistic level by carefully choosing the initial data and its approximation.
%%%%%
\subsection{Probabilistic well-posedness in the three dimensional subcritical case}\label{sec:3DPWP}
 It is shown in \cite{burq2013probabilistic} that for a specific random initial data together with a specific approximation of this initial data the problem  \eqref{eq:NLWaves} with $(\beta,d)=(1,3)$ is well-posed in the sense described below. Let us consider the initial data 
\begin{equation}\label{eq:FGF}
	\begin{cases}
		u_0^\omega(x)&= \displaystyle\sum_{n \in \mathbb{Z}^3}   \dfrac{g_\omega(n)}{\left \langle n\right\rangle^\alpha} e^{i  2\pi n \cdot x},\\
		\dot{u}_0^\omega(x)&=\displaystyle\sum_{n \in \mathbb{Z}^3}   \dfrac{h_\omega(n)}{\left\langle n\right\rangle^{\alpha-1}} e^{i  2\pi n \cdot x}.
	\end{cases}
\end{equation}
In \eqref{eq:FGF} $g(n)$ and $h(n)$ are independent and identically distributed complex random variables whose real and imaginary parts are independent and distributed according to the same law $\theta$ which verifies the sub-Gaussianity condition, 
\begin{equation}\label{eq:SubGauss}
	\exists\, c >0, \; \forall \gamma \in \R, \; \int_\R e^{\gamma x} d\theta(x)\leq e^{c\gamma^2}.
\end{equation} 
The sub-Gaussianity condition \eqref{eq:SubGauss} imposes that the tails of $\theta$ are decaying at least as fast as a Gaussian distribution. Additionally, $h$ and $g$ satisfy the hermitian symmetry $\overline{g_\omega(n)}= g_\omega(-n)$ and similarly for $h_\omega$.
Interestingly, the fractional Gaussian fields obtained with a Gaussian distribution of Fourier modes are ubiquitous in the modeling of three dimensional fluid turbulence \cite{krai,prr} where they appear as a Gaussian model of turbulent velocity fields.

 It has been shown in \cite{burq2013probabilistic} that the well-posedness of the problem depends on how one approximates the initial condition \eqref{eq:FGF}. Indeed, for $\alpha \in (3/2,2)$, one has $$ (u^\omega_0,\dot{u}_0^\omega) \in\Hs[1]{3} $$
almost surely for any $s \in (0,\alpha -3/2)$. One can therefore apply \cite[Theorem 1.28]{Tzvetkov2019random} and get 
the following statement. 
\begin{thm}[Pathological approximation]
	\label{thm:pathologicalApprox}
	Let $\alpha \in(3/2,2)$ and $s\in (0,\alpha-3/2)$. 
	For almost every $\omega$, there exists a sequence 
	$$ u^\omega_N \in C^0(\R, C^\infty(\T^3)), \quad N=1,2, \dots $$
	such that 
	$$(\partial^2_t -\Delta) u^\omega_N +(u_N^\omega)^3=0, $$
	with 
	$$ \lim_{N \rightarrow \infty } \| (u_N^\omega(0)-u_0^\omega, \partial_t u_N^\omega(0) -\dot{u}_0^\omega  )\|_{ \Hs[1]{3}}  =0$$
	but for all $T>0$,
	\[ \lim_{N \rightarrow \infty } \|( u_N^\omega(t), \partial_t u_N^\omega(t))   \|_{L^\infty\left( [0,T],  \Hs[1]{3}\right)}  =+\infty . \]
\end{thm}
Theorem \ref{thm:pathologicalApprox} does not rule out the existence of an approximation of the initial condition leading to convergence in $C^0([0,T], \Hs[1]{3})$ toward a unique and global weak solution of the non-linear wave equation and indeed, one has
\begin{thm}[Probabilistic well-posedness, \cite{burq2013probabilistic}]
	\label{thm:ProbaWP}
	Let  $\theta$ verify \eqref{eq:SubGauss}, $\alpha \in(3/2,2)$ and $s\in (0,\alpha-3/2)$. Set
	\begin{equation}\label{nachalna}
		\begin{array}{cr}
			u_{0,N}^\omega(x)= \displaystyle\sum_{ |n|\leq N}   \dfrac{g_\omega(n)}{\left\langle n\right\rangle^\alpha} e^{i  2\pi n \cdot x},&  v_{0,N}^\omega(x)=\displaystyle\sum_{ |n|\leq N}   \dfrac{h_\omega(n)}{\left\langle n\right\rangle^{\alpha-1}} e^{i  2\pi n \cdot x}.
		\end{array}
	\end{equation}
	Then there exists a set $\Sigma$ of probability $1$ such that for every $ \omega \in \Sigma$  the sequence $(u_N^\omega)_{N\geq 1}$ of solutions of  \eqref{eq:NLWaves} with smooth initial data given by \eqref{nachalna}
	converges when $N \rightarrow + \infty$ in $C^0(\R,H^s(\T^3))$ to a (unique) limit that satisfies \eqref{eq:NLWaves} in the sense of distributions.
\end{thm}
\begin{rk}\label{rk:RoleSubGauss}
	The sub-Gaussianity condition \eqref{eq:SubGauss} plays a crucial role in the proof of Theorem~\ref{thm:ProbaWP}. Loosely speaking, it ensures that the random series \eqref{eq:FGF} obeys large deviation estimates with Gaussian tails. Such estimates are what allows one to control the whole sequence of truncated solutions at once, and hence to reach a conclusion valid for almost every $\omega$; we refer to \cite{burq2013probabilistic} and \cite{Tzvetkov2019random} for the details. If $\theta$ has heavy tails, only finitely many moments of $g_\omega(n)$ and $h_\omega(n)$ are available, these estimates are no longer at our disposal, and such laws remain, for the time being, beyond the reach of our understanding.
\end{rk}
This is a limitation of the available proofs, and not an indication that the conclusion of Theorem~\ref{thm:ProbaWP} should fail for heavy tails distribution. One may therefore wonder whether \eqref{eq:SubGauss} can be relaxed. The numerical results presented in this article suggest that it can, and lead us to propose the following conjecture.
\begin{conj}\label{conj:HeavyTail}
	Theorem~\ref{thm:ProbaWP} remains valid when the sub-Gaussianity condition \eqref{eq:SubGauss} is replaced by the sole assumption that $\theta$ is centered with finite variance.
\end{conj}
We believe that Conjecture 1.5 is a challenging open problem. Its counterpart in the one dimensional fractional setting, stated as Conjecture~\ref{conj:HeavyTail1d} below, is the statement that our numerical experiments are designed to probe; the corresponding results are presented in Sections~\ref{sec:PWPNum} to \ref{sec:LWPNum}.

The stark contrast between the conclusions of Theorems \ref{thm:pathologicalApprox} and \ref{thm:ProbaWP} enlightens the very sensitive dependence on the choice of approximation of the initial conditions. 
We observed that for a large set of initial data in a measure theoretic sense  the approximation of initial data by Fourier truncation led to a converging sequence. On the other hand it turns out that there is another large set (in a topological sense) of initial data for which  the approximation of initial data by Fourier truncation leads to divergent sequences. More precisely, we have the following statement.
\begin{thm}
	[Pathological initial data  \protect{\cite{chenmin}}]
	Let $0<s<\frac12$. There is a dense set $S\subset H^s(\T^3) \times H^{s-1}(\T^3)$ such that for every $(f,g) \in S$, the sequence $(u_N)_{N\geq 1}$ of smooth solution of 
	$$  \partial_t^2u-\Delta u +u^3=0$$
	with (smooth) data 
	$$ 
	\begin{array}{cc}
		u_0^N(x) = \displaystyle  \sum_{|n| \leq N } \hat{f}(n) e^{2i\pi n \cdot x}, &\displaystyle  \dot{u}_0^N(x) =  \sum_{|n| \leq N } \hat{g}(n) e^{2i\pi n \cdot x},
	\end{array}
	$$
	does not converge. More precisely, for every $T>0$,
	$$  \lim_{N \rightarrow \infty} \| u_N(t)\|_{L^\infty([0,T], H^s(\T^3)    )} = + \infty.$$
\end{thm}
%
%All the existing results cited above concern the energy subcritical case. 
%
\subsubsection{The energy super-critical case}
Let us point out that all the results cited above concern the energy subcritical case. 
Little is known concerning the super-critical $\beta< \frac d4$ case. While it should be possible to show local well posedness in $\Hs{d}$ for $s>s_c$, the existence of finite time blowup of smooth solution is not excluded and remains an important open question. In this direction, it has recently been proven \cite{merleblowup} that the closely related super-critical defocusing nonlinear Schr\"odinger equation blows up in finite time. We therefore expect global existence for weak solutions only. We expect moreover that the ill-posedness results still holds and will provide numerical evidence that probabilistic well-posedness holds at least locally in time. 
\\[0.4cm]

For further references concerning probabilistic well-posedness, see \cite{Tzvetkov2019random} and \cite{Oh1,Oh2}. In the present article, we show that some of these fine properties of the non linear wave equation can be simulated numerically. For questions of numerical cost, we focus only on the one dimensional case. The study of a fractional dispersion allows us to work with a positive critical regularity $s_c=1/2-\beta$ provided that $\beta$ is small enough and therefore avoid renormalization. Introducing a fractional dispersion however makes the equation non-local and thus harder to study rigorously. The present numerical approach enables us to overcome this difficulty and to explore regimes that are inaccessible with the current analytical tools. Two such regimes are of interest here. The first one is the energy super-critical regime $\beta<\beta_c$, for which no global probabilistic well-posedness result is available. The second one, and the one we wish to put forward, is the regime of heavy tailed initial data, for which the sub-Gaussianity condition \eqref{eq:SubGauss} fails and which is the object of Conjecture~\ref{conj:HeavyTail}. 

\subsubsection{The one dimensional case}\label{sec:Transition1dNum}
Numerically, exploring one space dimension is more accessible than three. However, the previously exposed results for the three dimensional cubic case must be adapted. We expect the Theorems \ref{thm:pathologicalApprox} and \ref{thm:ProbaWP} to hold for a critical regularity $s_c=d/2-\beta$ meaning that in dimension $1$, we have to take $\beta <1/2$ for the critical regularity to remain positive. In addition, since the natural functional setup is $\Hs{ }$ we are going to study the well posedness of the non-linear wave equation \eqref{eq:NLWaves} with initial condition

\begin{equation}\label{eq:GoodInit}
	\begin{cases}
		u_0^\omega(x)&= \displaystyle\sum_{n \in \mathbb{Z}}   \dfrac{g_\omega(n)}{\left\langle n\right\rangle^\alpha} e^{i 2\pi n  x},\\
		\dot{u}_0^\omega(x)&=\displaystyle\sum_{n \in \mathbb{Z}}   \dfrac{h_\omega(n)}{\left\langle n\right\rangle^{\alpha-\beta}} e^{ i 2\pi n x},
	\end{cases}
\end{equation}
which belongs almost surely to $\Hs{}$ for any $0<s<\alpha -1/2 $ provided that $g(n)$ and $h(n)$ are centered i.i.d. random variables of finite variance. Note that this membership, which is the only place where the law $\theta$ enters the description of the initial data, requires nothing beyond a finite variance.

The statement that our numerical experiments are designed to test is the one dimensional fractional counterpart of Theorem~\ref{thm:ProbaWP} and of Conjecture~\ref{conj:HeavyTail}.
\begin{conj}\label{conj:HeavyTail1d}
	Let $d=1$, $0<\beta<1/2$ and let $\theta$ be centered with finite variance, not necessarily satisfying \eqref{eq:SubGauss}. Let $\alpha \in (1/2, 1-\beta)$, so that $\alpha-1/2 < 1/2-\beta$, and let $0<s<\alpha-1/2$. Then, for almost every $\omega$, the sequence $(u_N^\omega)_{N \geq 1}$ of solutions of \eqref{eq:NLWaves} with the Fourier truncated data \eqref{eq:WPInit} converges
	\begin{itemize}
		\item in $C\left(\R,\Hs{}\right)$ if $\beta>\beta_c=1/4$, that is in the energy sub-critical case;
		\item in $C\left([0,T],\Hs{}\right)$ for some $T=T(\omega)>0$ if $\beta<\beta_c=1/4$, that is in the energy super-critical case.
	\end{itemize}
\end{conj}
Conjecture~\ref{conj:HeavyTail1d} combines two extensions of Theorem~\ref{thm:ProbaWP}: the passage from the three dimensional cubic wave equation to the one dimensional fractional one, and the removal of the sub-Gaussianity condition \eqref{eq:SubGauss}. The distinction between the two items reflects the role of the conserved Hamiltonian \eqref{eq:HamiltonianNLWave}. In the energy sub-critical case $\beta>\beta_c$ the energy controls the norm in which the solution is measured, and one expects the local statement to be iterated globally in time, exactly as in Theorem~\ref{thm:ProbaWP}. In the energy super-critical case $\beta<\beta_c$ this control is lost and, as pointed out above, finite time blowup of smooth solutions is not excluded; we therefore only conjecture a local in time statement. Let us finally point out that our simulations, run up to the short time $t_s$, probe the local in time content of both items and say nothing directly about the global in time claim of the first one.

 In particular, we are going to observe that the choice of sequence approximating the initial condition \eqref{eq:GoodInit} plays a major role in the well posedness of the problem. 
Indeed, we expect that the above mentioned probabilistic well posedness result should also hold in one dimension if the initial data \eqref{eq:GoodInit} is approximated by
\begin{equation}\label{eq:WPInit}
	\begin{cases}
		&u_{0,N}^\omega(x)= \displaystyle\sum_{|n|\leq N}   \dfrac{g_\omega(n)}{\left\langle n\right\rangle^\alpha} e^{i 2\pi n  x},\\
		&\dot{u}_{0,N}^\omega(x)=\displaystyle\sum_{|n|\leq N}   \dfrac{h_\omega(n)}{\left\langle n\right\rangle^{\alpha-\beta}} e^{i  2\pi n x}.
	\end{cases}
\end{equation}
We will provide numerical evidence of norm inflation caused by a pathological approximation of  $(u_0^\omega,\dot{u}_0^\omega)$ by taking
\begin{equation}\label{eq:PathologicalInit}
	\begin{cases}
		&w_{0,N}^\omega(x)=  p_N^s(x) +\displaystyle\sum_{|n|\leq N}   \dfrac{g_\omega(n)}{\left\langle n\right\rangle^\alpha} e^{i 2\pi n  x},\\
	&	\dot{w}_{0,N}^\omega(x)=\displaystyle\sum_{|n|\leq N}   \dfrac{h_\omega(n)}{\left\langle n\right\rangle^{\alpha-\beta}} e^{i  2\pi n x}.
	\end{cases}
\end{equation}
Here $p_N^s$, is a perturbation going to zero in $H^s(\mathbb{T})$ while concentrating to one point and will be of the same form as the one used in the proof of ill-posedness in \cite{Tzvetkov2019random}. Typically, we will consider
\begin{equation}\label{eq:Pertubation}
	p_N^s(x)=  \begin{cases}
		&c\dfrac{N^{1/2-s}}{\log(1+N) }\exp\left(1- \dfrac{1}{1-4N^2(x-\frac12)^2} \right), \; \text{if}\;  \left|x-\frac12 \right|<\frac{1}{2N},\\
		&0,\; \text{else}.
	\end{cases}
\end{equation}
The perturbation is a bump localized around $x=1/2$, of typical width $1/N$ and height $c N^{1/2-s}/\log (1+N) \to \infty $ when $ N\to \infty$ and $s<1/2$. While the perturbation concentrates in a point, it is vanishing in $H^s(\T)$ since
$$ \|p_N^s\|_{H^s(\T)} \lesssim \frac{1}{\log(1+ N)}. $$
This ensure that both \eqref{eq:WPInit} and \eqref{eq:PathologicalInit} converge toward the same initial condition \eqref{eq:GoodInit} in $\Hs{}$ for any $s< \alpha -1/2$. 

In addition of studying the energy subcritical and super critical regimes, we will also explore two different distributions for $g(n)$ and $h(n)$, one verifying the sub-Gaussianity assumption \eqref{eq:SubGauss} and another one violating it. In particular, we will present results for a Gaussian and a Student distribution,
\begin{equation}\label{eq:DistribFourier}
	d\theta_g(x)=  \dfrac{1}{\sqrt{2\pi}} e^{-\frac{x^2}{2}} dx, \qquad d\theta_s(x)=   \dfrac{2}{\pi} \left( 1+x^2\right)^{- 2} dx.
\end{equation}
The law $\theta_s$ is the Student distribution with three degrees of freedom, rescaled so as to have unit variance. Both laws in \eqref{eq:DistribFourier} are centered and of unit variance, so that in both cases the initial condition \eqref{eq:GoodInit} almost surely lies in $\Hs{}$ for any $0<s<\alpha-1/2$, with the same second moment
$$ \mathbb{E}\Big[ \|u_0^\omega\|^2_{H^{s}(\T)}\Big] = 2\sum_{n \in \mathbb{Z}} \left\langle n \right\rangle^{2s-2\alpha} <+\infty ,$$
and similarly for $\dot{u}_0^\omega$ in $H^{s-\beta}(\T)$. The Gaussian law $\theta_g$ satisfies the assumption \eqref{eq:SubGauss} whereas $\theta_s$ clearly violates it: its tails decay only like $|x|^{-3}$, so that $\int_\R |x|^p \,d\theta_s(x)$ is finite if and only if $p<3$. In particular, in the Student case the random variable $\|u_0^\omega\|^2_{H^{s}(\T)}$ has a finite expectation but an infinite variance, while all its moments are finite in the Gaussian case. The two settings are thus markedly different from a probabilistic point of view, even though the initial data enjoy the same almost sure Sobolev regularity.

\subsection{Organization of the numerical investigation}\label{sec:OrgaPaper}

Our numerical investigation is threefold. First of all, in Section~\ref{sec:PWPNum} we will present numerical results concerning the probabilistic well posedness of \eqref{eq:NLWaves} in both energy subcritical $(\beta >1/4)$ and supercritical $(\beta <1/4)$ regimes and for both Gaussian and Student distribution of Fourier modes \eqref{eq:DistribFourier}. The initial data \eqref{eq:GoodInit} will therefore be below the critical regularity for deterministic well posedness ($\alpha -1/2 < 1/2-\beta$) and approximated by Fourier truncation \eqref{eq:WPInit}.

 Then in Section~\ref{sec:NormInflationNum}, comparing with the case of the probabilistically well posed setting, we will study the norm inflation phenomenon in both energy subcritical $(\beta >1/4)$ and supercritical $(\beta <1/4)$ regimes and for both Gaussian and Student distribution of Fourier modes \eqref{eq:DistribFourier}. This also corresponds to the low regularity regime ($\alpha -1/2 < 1/2-\beta$) but this time with the pathological approximation \eqref{eq:PathologicalInit} of the initial condition \eqref{eq:GoodInit}.

Finally, in Section~\ref{sec:LWPNum}, we will study the deterministic well posedness regime ($\alpha -1/2 > 1/2-\beta$) for which the regularity of the initial condition \eqref{eq:GoodInit} is high enough to ensure deterministic local well-posedness of the considered Cauchy problem. In this case, both approximations \eqref{eq:WPInit} and \eqref{eq:PathologicalInit} should lead to the same well behaved solution of \eqref{eq:NLWaves} in the large $N$ limit. We will thus test this for both Gaussian and Student distribution of Fourier modes \eqref{eq:DistribFourier}. Verifying this numerically is therefore a sanity check regarding the consistency of our numerical simulations.

In each of these three sections the two laws \eqref{eq:DistribFourier} are treated on the same footing, and the comparison between them is the main outcome of this work: in all three regimes, and in both the energy sub-critical and the energy super-critical case, the heavy tailed simulations reproduce the sub-Gaussian ones without qualitative change. This is the numerical evidence we can offer in favour of Conjecture~\ref{conj:HeavyTail1d}, and indirectly of Conjecture~\ref{conj:HeavyTail}.

\subsubsection{Numerical setting for probabilistic well-posedness}
In Section \ref{sec:PWPNum} we will present results concerning the probabilistic well posedness of \eqref{eq:NLWaves} advocating in favor of an extension of Theorem~\ref{thm:ProbaWP} to our one dimensional fractional setting in both energy sub-critical and energy super-critical regimes. In order to do so, we run numerical simulations of \eqref{eq:NLWaves} in one dimension with a given realization of the random initial condition \eqref{eq:GoodInit}, for each of the two laws \eqref{eq:DistribFourier}, with low regularity ($\alpha -1/2< 1/2-\beta $) and approximated by \eqref{eq:WPInit}. As exposed in Section \ref{sec:PWPNum}, for both $\beta=1/3> \beta_c =1/4$ (energy sub-critical) and $\beta=1/8 \leq \beta_c$ (energy super-critical), the numerical approximation of the solution remains bounded and converges in $C([0,T],\Hs{})$ in the limit where the approximation \eqref{eq:WPInit} converges to \eqref{eq:GoodInit}. At the level of the initial condition, the speed of convergence of the approximation in $L^2(\Omega,\Hs{})$ for any $0<s<\alpha -1/2$ is given by
\begin{equation} \label{eq:ConvPWPInitCond}
	\| (u_0^\omega-u_{0,N}^\omega , \dot{u}_0^\omega-\dot{u}_{0,N}^\omega  ) \|_{L^2_\omega \mathcal{H}_{s,\beta}}  \leq C  N^{s-(\alpha-\frac12)},
\end{equation}
where $C$ is a constant depending on $s$ and $\alpha$ only. The estimate \eqref{eq:ConvPWPInitCond} has a practical implication concerning the numerical simulations. It indicates that the closer $s$ is to zero, the faster the approximation converges in $L^2(\Omega,\Hs{})$ which is a remark of importance when choosing the parameters for numerical simulations. As a trade-off between convergence speed and strength of the studied norm, we will take $ s= (\alpha-1/2)/3$ and $\alpha=0.6$ in the numerical simulations concerning probabilistic well-posedness.

\subsubsection{Numerical setting for norm inflation}
 Comparing with the case of the probabilistically well posed setting, we will provide numerical evidence in Section \ref{sec:NormInflationNum} that the solution of \eqref{eq:NLWaves} with pathological approximation \eqref{eq:PathologicalInit} of the initial condition \eqref{eq:GoodInit} gets arbitrarily large in $C([0,T],\Hs{})$ for any $T>0$ when $N$ tends to infinity. To this end we note for comparison that for the pathological approximation \eqref{eq:PathologicalInit} we have
\begin{equation} \label{eq:ConvPathologicalInitCond}
	\| (u_0^\omega-w_{0,N}^\omega , \dot{u}_0^\omega-\dot{w}_{0,N}^\omega  ) \|_{L^2_\omega\mathcal{H}_{s,\beta}} \leq C N^{s-(\alpha-\frac12)} +K (\log N )^{-1} ,
\end{equation}
where $C$ is the same constant as in \eqref{eq:ConvPWPInitCond} and $K$ depends only on the shape of the perturbation. The pathological approximation \eqref{eq:PathologicalInit} is thus converging to \eqref{eq:GoodInit} only logarithmically. In order to compare the behavior of the pathological approximation \eqref{eq:PathologicalInit} with the Fourier truncation approximation \eqref{eq:WPInit}, we will also measure the $\Hs{}$ norm with $ s= (\alpha-1/2)/3$ and $\alpha=0.6$ when investigating norm inflation.

\subsubsection{Numerical setting for deterministic well-posedness}
As sketched in Section \ref{sec:LocalWP} we expect deterministic well posedness of the problem \eqref{eq:NLWaves} for high enough Sobolev regularity of the initial data. More precisely, we expect local well posedness whenever the initial condition $(u_0, \dot{u}_0)$ (not necessarily random) lies in $\mathcal H_{\sigma,\beta}(\T)$ for some $\sigma >1/2-\beta$. For deterministic well posedness to hold, the dependence of the solution on the initial condition must be continuous in the same Sobolev space $\mathcal H_{\sigma,\beta}(\T)$. In our case of random initial data \eqref{eq:GoodInit}, we recall that for any $0< \sigma < \alpha -1/2$ the initial condition $(u_0^\omega, \dot{u}_0^\omega)$ lies almost surely in $\mathcal H_{\sigma,\beta}(\T)$. This implies that as long as $1/2-\beta < \alpha -1/2$, one can find $\sigma >1/2-\beta$ such that $(u_0^\omega, \dot{u}_0^\omega) \in \mathcal H_{\sigma,\beta}(\T)$ almost surely. Then, because of the choice of perturbation $p_N^s$, the approximation \eqref{eq:PathologicalInit} is converging in $\mathcal H_{\sigma,\beta}(\T)$  toward the initial data \eqref{eq:GoodInit} provided that $\sigma \leq s$. We have indeed $\|p_N^s\|_{H^\sigma} \simeq N^{\sigma-s}/\log N$. For this convergence in $\mathcal H_{\sigma,\beta}(\T)$ to hold for some $\sigma > 1/2-\beta$, one need to take $1/2-\beta< s$. 

Therefore, provided that the initial condition \eqref{eq:GoodInit} enjoys enough regularity (that is $1/2 -\beta < \alpha -1/2$) and that the perturbation $p_N^s$ vanishes in $H^\sigma(\T)$ for some $\sigma$ such that $1/2-\beta< \sigma < \alpha- 1/2$  (that is $1/2-\beta<s $), there exists $\sigma$ verifying $1/2-\beta <\sigma \leq s< \alpha-1/2 $ such that the initial data \eqref{eq:GoodInit} lies in $\mathcal H_{\sigma,\beta}(\T)$ and the pathological approximation \eqref{eq:PathologicalInit} converges to \eqref{eq:GoodInit} in $\mathcal H_{\sigma,\beta}(\T)$. In the numerical study of deterministic well posedness, we will place ourselves at the endpoint $\sigma = s$ and thus study deterministic well posedness in $\Hs{}$ with the constraint $1/2 -\beta < s < \alpha -1/2$.  

The purpose of Section \ref{sec:LWPNum} is to present numerical results showing that \eqref{eq:NLWaves} is locally well posed in $\Hs{}$ for both approximations \eqref{eq:PathologicalInit} \& \eqref{eq:WPInit} of the initial data \eqref{eq:GoodInit}  whenever $1/2-\beta <s< \alpha -1/2$. A lack of convergence in this regime would signal an inappropriate numerical integration of the dynamics. This deterministic well posedness regime thus constitutes a fail-safe for our numerical simulations.   
In the numerical simulations we will consider $s=\gamma (\alpha-1/2)$. The constraint $1/2-\beta <s< \alpha -1/2$ can thus be transposed to the value $\gamma$. In particular, in the study of deterministic well-posedness we will take 
\begin{equation} \label{eq:ConstraintGamma}
	\dfrac{\frac12-\beta }{ \alpha - \frac12}< \gamma= \frac{1}{2}\left(1 + \dfrac{\frac12-\beta }{ \alpha - \frac12}\right) <1,
\end{equation}
which is the midpoint between $(\frac12-\beta )/( \alpha - \frac12)$ and $1$. In addition, we will consider $\alpha= 0.98$ and $\beta= 1/3$ for the energy subcritical case and $\beta =1/8$ for the energy super critical one.

\begin{figure}[htb]
	\centering
	\includegraphics[width=0.7\linewidth]{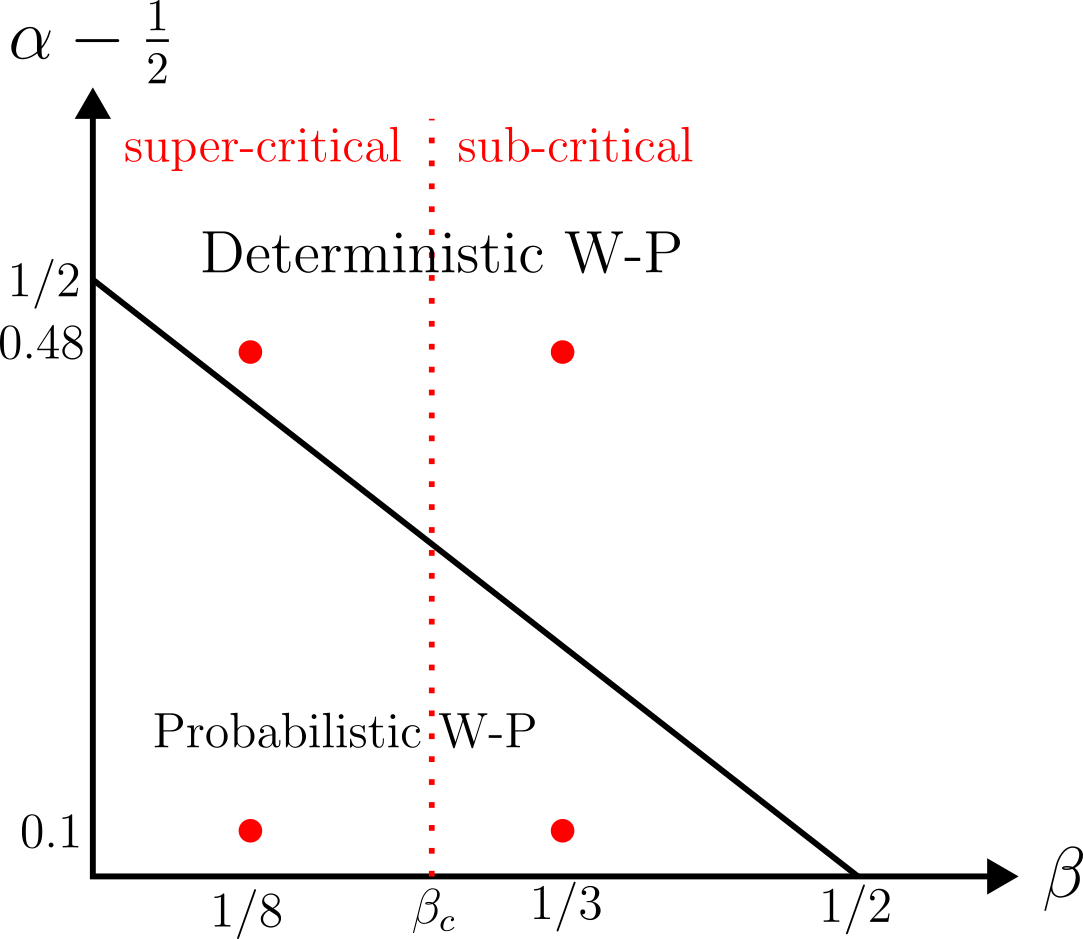} 
	\caption{Illustration of different expected behaviors of \eqref{eq:NLWaves} when the Sobolev regularity of the initial data and dispersion vary. The red dotted line marks the distinction between energy super-critical $\beta< \beta_c =1/4$ and sub-critical regions. The solid black line $ \alpha-1/2=1/2-\beta$ is the boundary between deterministic well posedness and probabilistic well-posedness. The red dots correspond to the numerical experiments that we will present in the subsequent sections }
	\label{fig:Regimes}
\end{figure}
For the sake of clarity, we summarize the different behaviors of the solution of \eqref{eq:NLWaves} in Figure~\ref{fig:Regimes}, describing when we expect deterministic or probabilistic well posedness to hold and the frontier between energy sub-critical and energy super-critical regimes symbolized by a red dashed line. We also place red dots in the $(\alpha,\beta)$ plane corresponding the numerical experiments that we will present in sections \ref{sec:PWPNum}, \ref{sec:NormInflationNum} and \ref{sec:LWPNum}. As sketched in Figure~\ref{fig:Regimes}, the phenomena we seek to illustrate are completely determined by the parameters $\alpha$ (measuring the roughness of the initial data) and $\beta$ (measuring the strength of dispersion). The Sobolev spaces we are working with only appear as technical tools.

Finally, in Section~\ref{sec:HamilCons} we will provide numerical results on the conservation of the Hamiltonian \eqref{eq:HamiltonianNLWave} in order to probe the precision of our numerical procedure. Note that the Hamiltonian \eqref{eq:HamiltonianNLWave} might diverge in the large $N$ limit, we therefore compute the error on the conservation of the Hamiltonian relatively to its value at initia time.

\section{Numerical investigation}\label{sec:Num}
\subsection{Numerical methods}\label{sec:NumMethods}
In order to capture the instantaneous growth of Sobolev norms expected with norm inflation, we need a robust numerical method together with a high resolution of numerical simulations. In order to account for the loss of regularity in the large $N$ limit, the time step used in numerical simulations will depend on $N$ and in particular vanish as $N$ gets infinitely large. The spatial discretization will also obviously depend on $N$. In the sequel we will use a so-called filtered trigonometric integrator.

\subsubsection{Time integration}\label{sec:TrigoScheme}
Trigonometric integrators are a class of numerical integration schemes known for their good conservation of the Hamiltonian structure due to their symplectic nature \cite{cohen2008conservation}. In the following, the time interval $[0,t_s]$ will be discretized with a time step $\tau$ and for $p\in \mathbb{N}$ we will denote by $(u_p, \dot{u}_p)$ the numerical approximation of $(u,\partial_t u)$ at time $t_p=p\tau$.

We are going to use the trigonometric integrator with filtered non-linearity proposed in \cite{gauckler2019trigonometric}. For convenience, let us write the non-linear term as $f(u)=-u^3$ and introduce the filtered non-linearity
$$ \tilde{f}(u)= \sinc^2(\tau \D{\beta}) f\left(  \sinc(\tau \D{\beta})   u \right) ,$$
where $ \sinc(\tau \D{\beta}) $ is understood as the Fourier multiplier $ \sinc(\tau |2\pi k|^{\beta}) $ and $\sinc(x) = \sin(x)/x$. Then, the trigonometric integration scheme writes,
\begin{equation}\label{eq:NumScheme}
	\begin{cases}
		u_{p+1}&=\cos(\tau \D{\beta})u_p+\tau \mathrm{sinc}(\tau \D{\beta}) \dot{u}_p +\frac 12 \tau^2\mathrm{sinc}(\tau \D{\beta}) \tilde{f}(u_p)   ,\\
		\dot{u}_{p+1}&=-\D{\beta} \sin(\tau \D{\beta}) u_p+ \cos(\tau \D{\beta}) \dot{u}_p + \frac 12 \tau \left(  \cos(\tau \D{\beta} )  \tilde{f}(u_p) +  \tilde{f}(u_{p+1})  \right).
	\end{cases}
\end{equation}
For any finite value of $N$, the initial datum \eqref{eq:WPInit} and \eqref{eq:PathologicalInit} are smooth, ensuring convergence of the numerical scheme. In order to have a sound numerical approximation, it is necessary to have a space discretization in accordance with the number of energy containing Fourier modes. The integration scheme \eqref{eq:NumScheme} is iterated $P$ times, meaning that the solution is approximated for times between $0$ and $t_s= P \tau$. Since norm inflation is instantaneous in the limit, we don't need to explore large simulation times and will typically take $t_s=10^{-2}$.

\subsubsection{Space discretization} \label{sec:SpaceDiscr}
Given the periodic nature of the problem and the efficient computation of the linear operators in Fourier space using fast Fourier transform, we adopt a pseudo-spectral method. We discretize the torus $\T= [0,1]$ using $M$ collocation points, meaning that the spatial step size is $h=1/M$ and Fourier modes take discrete values such that $-M/2+1\leq n\leq M/2$. The value of $M$ must be large enough to resolve both the truncated fractional Gaussian field and the perturbation $p_N^s$. In practice in our simulations we go up to $M= 2^{25}$ allowing to reach a truncation level $N= 2^{23}$. In addition, the cubic nonlinearity is fully dealiased. Dealiasing is performed by zero padding on a grid of size $2M$ because of the cubic nature of the nonlinearity. 

\subsection{Observables of interest}\label{sec:Observables}
The precision of numerical integration is crucial to study such fine behaviors of the nonlinear wave equation \eqref{eq:NLWaves}. We will thus keep track of the conservation of the Hamiltonian \eqref{eq:HamiltonianNLWave} by the numerical flow as a substitute for precision of the numerical integration. The Hamiltonian \eqref{eq:HamiltonianNLWave} will be approximated numerically by
\begin{equation}\label{eq:DiscrHamil}
	\mathcal{H}_{\text{num}}[u,\dot{u}](p\tau)= \sum_{ -\frac M2 < n \leq \frac M2} |\widehat{\dot{u}}_p(n)|^2+ |2 \pi n |^{2\beta} |\widehat{u}_p(n)|^2+ \dfrac{1}{2M} \sum_{0 \leq q <  M  } u_p(qh )^4,
\end{equation}
The first two terms are the numerical approximations of the square $L^2$ norms of $\partial_t u $ and $D_x^\beta u$ based on Parseval's theorem and the last term is the Riemann sum approximation of the integral of $u^4/2$.  
In the following numerical study, the most important part being the Fourier truncation level $N$ we will keep implicit the dependence on $p$ and $M$ in order to avoid heavy notations and refer to the numerical approximation of the solution as $(u_N(t),\dot{u}_N(t))$ for $t \in \tau \lbrace 0, \dots, P\rbrace$. For any finite value of $N$, the flow of \eqref{eq:NLWaves} is smooth and thus conserves the Hamiltonian \eqref{eq:HamiltonianNLWave}. Here, we will track conservation of the discretized Hamiltonian \eqref{eq:DiscrHamil} by the numerical integration, such conservation should hold exactly in the limit $M, \tau^{-1}\to \infty$ for any finite $N$. In particular, since the Hamiltonian \eqref{eq:HamiltonianNLWave} need not be finite in the limit $N\to \infty$ we will track the relative error on the discretized Hamiltonian defined by 
\begin{equation}\label{eq:NumError}
	e_\infty^N= \max_{ p\in[0,P]} \left| \dfrac{\mathcal{H}_{\text{num}}[u_N,\dot{u}_N](\tau p)}{\mathcal{H}_{\text{num}}[u_N,\dot{u}_N](0)}-1 \right|.
\end{equation}
In addition of tracking the evolution of the error \eqref{eq:NumError} with $N$, we want to investigate the behavior of the solution in the Sobolev space $\Hs{}$. In the next sections we are therefore going to compute the Sobolev norms
\begin{equation}
	\begin{cases}
		&\SN=\displaystyle \sqrt{ \sum_{-\frac M2 < n \leq \frac M2} \left\langle n\right\rangle^{2s-2\beta} |\widehat{\dot{u}}_N(t,n)|^2  }+ \sqrt{ \sum_{-\frac M2 < n \leq \frac M2}\left\langle n\right\rangle ^{2s} |\widehat{u}_N(t,n)|^2   } \\
		& \phantom{\SN} \simeq \|(u_N^\omega(t),\dot{u}_N^\omega(t))\|_{ \Hs{}}, \\
		& \Sinf= \underset{ t \in[0,t_s]}{\max} \, \SN  \simeq \|(u_N^\omega,\dot{u}_N^\omega)\|_{L^\infty([0, t_s], \Hs{})}.
	\end{cases}\label{eq:SobolNum}
\end{equation}
In order to make a distinction between the two different approximations \eqref{eq:WPInit} and \eqref{eq:PathologicalInit} of the initial condition , we will change the notation for the Sobolev norms in the pathological approximation case  \eqref{eq:PathologicalInit} to be $\NN$ and $\Ninf$ respectively.
The numerical computation of the norms $\Sinf$ and $\Ninf$ defined in \eqref{eq:SobolNum} allows us to address the boundedness of the solution in $C([0,t_s],\Hs{})$ in the large $N$ limit but is not sufficient to decide whether or not the sequence of solutions $(u_N,\dot{u}_N)$ converges in $C([0,t_s],\Hs{})$. To this end, we will also probe the convergence to zero of the sequence $(u_N-u_{N/2})_{N \in \mathbb{N}}$ by computing 
\begin{equation*}
	\Delta \Sinf = \underset{ t \in[0,t_s]}{\max} \,  \left(\displaystyle \sqrt{ \sum_{-\frac M2 < n \leq \frac M2} \left\langle n \right\rangle^{2s-2\beta} \big|\widehat{\dot{u}}_N(t,n) -\widehat{\dot{u}}_{\frac N2}(t,n)\big|^2  } \right. \\
	\left. + \sqrt{ \sum_{-\frac M2 < n \leq \frac M2}  \left\langle n \right\rangle^{2s} \big|\widehat{u}_N(t,n) -\widehat{u}_{\frac N2}(t,n)\big|^2   }\right) ,
\end{equation*}
which is a discrete approximation of $ \|(u_N^\omega-u_{ N/2}^\omega,\dot{u}_N^\omega-\dot{u}_{ N/2}^\omega)\|_{L^\infty([0, t_s], \Hs{})}$. Similarly, we will denote $\Delta \Ninf$ the same norm in the case of the pathological approximation \eqref{eq:PathologicalInit}. The convergence to zero of $ \Delta \Sinf$ or $ \Delta \Ninf$ alone is not sufficient to claim that $(u_N)_{N\in \mathbb{N}}$ or $(w_N)_{N\in \mathbb{N}}$ are Cauchy in $C([0, t_s], \Hs{})$. However, such a behavior goes in the direction of convergence in $C([0, t_s], \Hs{})$ while a failure of $\Delta \Sinf$ or $\Delta \Ninf$ to converge to $0$ is a clear evidence of absence of convergence in $C([0, t_s], \Hs{})$.

For the reasons detailed in Section \ref{sec:OrgaPaper} the probed Sobolev norms will be of index  $s= \gamma ( \alpha -1/2)$ with $\gamma=1/3$ if $\alpha-1/2 < 1/2- \beta $ (probabilistic well posedness regime) and $\gamma = 1/2 (1+(1/2-\beta)/(\alpha-1/2)   )$ otherwise (deterministic well posedness regime).

\subsection{Choices of parameters for the simulations}\label{sec:ParamsSimuls}
As stated in Section~\ref{sec:OrgaPaper} we wish to illustrate the three following behaviors of \eqref{eq:NLWaves}: Probabilistic well posedness, norm inflation and deterministic well posedness. In each case we explore both the energy sub-critical and energy super-critical regimes together with Gaussian and Student distributions of Fourier modes \eqref{eq:DistribFourier}. Probabilistic well posedness and norm inflation will be investigated using $\alpha=0.6$ and $\beta =1/3>  \beta_c =1/4$ for the energy subcritical regime and $\beta =1/8<\beta_c$ for the energy super-critical one. For the two different values of $\beta$, we indeed have $\alpha-1/2 < 1/2-\beta$ which is the range where deterministic well-posedness breaks down. In the case of local well-posedness, we consider $\alpha=0.98$ and the same values of $\beta$ as for the previous cases.  
We explore values of $N=2^k$ with $k\in\lbrace 3, \dots, 23\rbrace$, for reasons of computation time we set the spatial resolution $M$ to be $2^{23}$ whenever $k\in\lbrace 4, \dots, 20\rbrace$ and $M=2^{25}$ for $k\in\lbrace 21, 22, 23\rbrace$. For each value of $N$ we integrate \eqref{eq:NLWaves} using the pseudo-spectral trigonometric integrator described in Section~\ref{sec:TrigoScheme} from $t=0$ up to $t_s= 10^{-2}$ with a time step $\tau_N$ depending on $N$. The time step $\tau_N$ should be small enough to resolve both dispersive and nonlinear effects. The dispersive time scale for the largest Fourier mode is $\tau_d = M^{-\beta}$. The nonlinear time scale $\tau_{\text{NL}} $ is inversely proportional to the amplitude of the solution. Estimating roughly the amplitude of the initial condition one gets a contribution from the perturbation $p_N^s$ and a contribution from the random part. The amplitude of $p_N^s$ scales as $N^{1/2-s}/\log N$ and the amplitude of the random part is estimated as the standard deviation of \eqref{eq:GoodInit} which is close to $\sqrt{\zeta(2\alpha)}$, where $\zeta$ is the Riemann $\zeta$ function. Since $\theta_g$ and $\theta_s$ have the same variance, this estimate, and therefore the time step, is the same for the two laws \eqref{eq:DistribFourier}. Therefore, at initial time, the nonlinear time scale is defined as $(\tau^0_{\text{NL}})^{-1} = \sqrt{\zeta(2\alpha)}+ N^{1/2-s}/\log N$. However, since the solution is quickly amplified in the case of norm inflation, a non linear time scale based on the initial condition shall not be small enough to resolve properly the time evolution. As a rule of thumb empirically validated by the monitoring of energy conservation presented in Section~\ref{sec:HamilCons}, we thus consider
\begin{equation}\label{eq:TimeStep}
	\tau_N = \min\left( \dfrac{\tau_d}{5}, \frac{(\tau^0_{\text{NL}})^{3}}{2}\right),
\end{equation} 
where the cubic dependence on $\tau^0_{\text{NL}}$ accounts for the amplification of the solution and will turn out to be a rather conservative choice.
Note that the time step goes to zero when $\alpha$ goes to $1/2$ or $N$ goes to infinity, this accounts for the low regularity of the limiting solution.
As we will also detail in Section~\ref{sec:HamilCons}, we perform numerical experiments with a refined time step $\tau/2$ together with a refined spatial discretization $2M$ in order to show that our results are unchanged upon this refinement.% We also performed simulations with a St\"ormer-Verlet integration scheme (data not shown) instead of \eqref{eq:NumScheme} without qualitative change in the presented results. 

 Importantly, probabilistic well-posedness holds for almost every realization $\omega$ of the initial data \eqref{eq:GoodInit}. Performing numerical simulations with a different realization of the initial condition will yield quantitatively different results than those exposed here but the conclusion will remain unchanged. In this regard, our numerical results are reproducible up to the choice of seed for the random number generator. Also, for each of the two laws \eqref{eq:DistribFourier}, a single realization of $(g_\omega,h_\omega)$ entering \eqref{eq:GoodInit} is drawn once and shared by all the simulations reported below, meaning that only $\alpha,\, \beta, \, N$ and $s$ change between different simulations performed with a given distribution. The monitoring of the Hamiltonian presented in Section~\ref{sec:HamilCons} shows that the choice \eqref{eq:TimeStep} of time step remains appropriate for the heavy tailed law $\theta_s$.

\subsection{Probabilistic well-posedness  }\label{sec:PWPNum}
\begin{figure}[htb]
	\centering
	\includegraphics[width=0.80\linewidth]{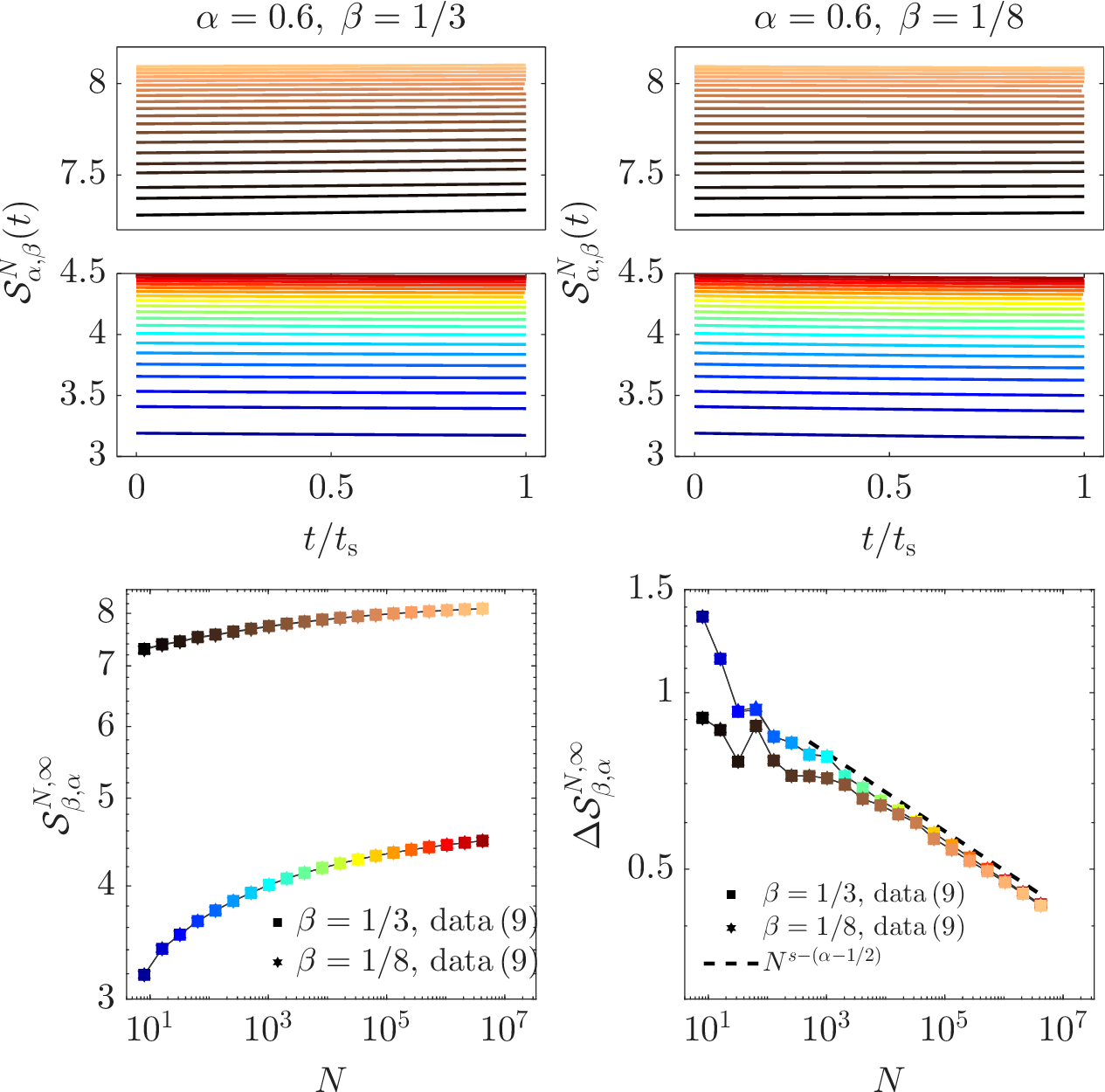}
	\caption{Probabilistic well posedness regime ($\alpha=0.6$), initial data approximated by the Fourier truncation \eqref{eq:WPInit}. The upper row shows the time evolution of $\SN$ defined by \eqref{eq:SobolNum} for different values of $N\in\lbrace 2^3, \dots 2^{23} \rbrace$ corresponding to the color of each line, matching the color of the markers in the bottom row. The left column of the first row is dedicated to the $(\alpha,\beta)=(0.6,1/3)$ energy subcritical case while the right one is dedicated to the energy super-critical $(\alpha,\beta)=(0.6,1/8)$ case. Within each column, the upper sub-panel corresponds to the Student law $\theta_s$, with a colormap going from black (smallest $N$) to light brown (largest $N$), and the lower sub-panel to the Gaussian law $\theta_g$, with a colormap going from blue (smallest $N$) to red (largest $N$); the same color convention is used in Figures~\ref{fig:NormInflation}, \ref{fig:LWP} and \ref{fig:ErrorHamil}. On the left panel of the bottom row, we show the evolution with $N$ of $\Sinf$ for $\beta=1/3$ with square markers and $\beta=1/8$ with star markers, each marker being one simulation; the upper group of markers corresponds to the Student law and the lower one to the Gaussian law. Since the $\Hs{}$ norm of the initial condition is independent of $\beta$, the cases $\beta=1/3$ and $\beta=1/8$ are indistinguishable by eye. Finally, the right panel of the bottom row shows the evolution of $\Delta\Sinf$ for $\beta=1/8$ (stars) and $\beta=1/3$ (squares) and for the two laws, we also superimpose the power law $N^{s-(\alpha-\frac12)}$ (black dashed line) in order to compare to the convergence of the approximation at initial time given by \eqref{eq:ConvPWPInitCond}. Once again, the cases $\beta=1/3$ and $1/8$ cannot be distinguished by eye, and so are the Gaussian and Student cases.   }
	\label{fig:PWP}
\end{figure}

We can now present the numerical illustration of probabilistic well posedness. The results of numerical simulations are displayed in Figure~\ref{fig:PWP}. They are obtained by performing simulations of \eqref{eq:NLWaves} with initial data \eqref{eq:WPInit} for increasing values of the truncation level $N$ and for the two laws \eqref{eq:DistribFourier} of the Fourier coefficients. Throughout Sections~\ref{sec:PWPNum} to \ref{sec:HamilCons} we use the following color code: the Gaussian case $\theta_g$ is displayed with a colormap ranging from blue (smallest $N$) to red (largest $N$), while the Student case $\theta_s$ is displayed with a colormap ranging from black (smallest $N$) to light brown (largest $N$). For each value of the truncation number $N$, we compute the Sobolev norm $\SN$ of the numerical solution and represent it as a function of time and for $\alpha=0.6$ in the energy sub-critical $\beta=1/3$ (upper left panels) and energy super-critical $\beta=1/8$ (upper right panels) cases; within each column, the upper sub-panel is devoted to the Student law and the lower one to the Gaussian law, following the color code detailed above. Then, we remove time out of the picture by computing for each value of $N$ the supremum of $\SN$, noted $\Sinf$ (bottom left panel) in order to  determine whether or not the solution is bounded in $C([0,T],\Hs{})$. Finally we probe the convergence in $C([0,T],\Hs{})$. The bottom right panel displays the behavior of $\Delta \Sinf$ which tends to show that $(u_N)_{N\in \mathbb{N}}$ is Cauchy and thus converges in $C([0,T],\Hs{})$. In the two panels of the bottom row, the results obtained with the two laws \eqref{eq:DistribFourier} are superimposed on the same axes. 

From the upper row of Figure~\ref{fig:PWP}, it is clear that for any fixed $N$ and for both laws \eqref{eq:DistribFourier}, the Sobolev norm $\SN$ has a very mild dependence on time and its initial value is independent of $\beta$ as expected from the choice of initial data. More importantly the four sub-panels of the upper row advocate for a uniform in time convergence of $\SN$ when $N$ increases (from blue to red in the Gaussian case, from black to light brown in the Student case). This uniform in time convergence is a feature of both energy sub-critical and super-critical \eqref{eq:NLWaves} and implies boundedness of the sequence $(u_N)_{N\in \mathbb{N}}$ in $C([0,t_s], \Hs{})$. More quantitatively, the left panel of the lower row of Figure~\ref{fig:PWP} shows the behavior of $\Sinf$ with $N$. Each marker corresponds to a simulation and the color of the markers depends on the Fourier truncation level $N$ and matches the colors of solid lines in the upper row. Probabilistic well posedness of the Cauchy problem \eqref{eq:NLWaves}\&\eqref{eq:GoodInit} can hold only if $\Sinf$ has a finite limit when $N$ goes to infinity. Even though the convergence with $N$ of $\Sinf$ might yet seem unclear from Figure~\ref{fig:PWP}, we stress that the $x$ and $y$ axes of the figure are logarithmically scaled, meaning that if $\Sinf$ diverges with $N$, it does slower than any power of $N$ since the data shows no linear behavior in log-log scale. In addition when examining the same data with logarithmic scale on the $x-$axis only, no linear behavior in terms of $\log N$ are found, meaning that if $\Sinf$ diverges, it does so slower than $\log N$. From a numerical perspective, a growth slower than a logarithm is hard to distinguish from a constant on the range of $N$ we can explore. In fact, as supported by the upper row of Figure~\ref{fig:PWP}, the slow convergence of $\Sinf$ can be accounted for by the slow convergence of the approximation of the initial condition in $\Hs{}$. Indeed, we remind that the approximation \eqref{eq:WPInit} converges at a rate $N^{s -(\alpha-1/2)}$ toward \eqref{eq:GoodInit}. With our parameters this gives a convergence rate of $ N^{-\mu}$ with $\mu \simeq 0.06 $ which is small and explains the slow convergence of $\Sinf$. Finally, in the right panel of the bottom row of Figure~\ref{fig:PWP}, we display the behavior of $\Delta \Sinf$ in order to decide if the sequence $(u_N)_{N \in \mathbb{N}}$ converges in $C([0,T],\Hs{})$. As shown in Figure~\ref{fig:PWP}, $\Delta \Sinf$ shows a convergence to zero as a negative power of $N$ corresponding to the convergence rate of the approximation of the initial condition. Even though the convergence of $\Delta \Sinf$ to zero does not necessarily imply that $(u_N)_{N \in \mathbb{N}}$  is a Cauchy sequence, this result goes toward such a property and thus convergence in $C([0,T],\Hs{})$ as expected for probabilistic well-posedness. Let us finally mention that our simulations only investigate the local in time probabilistic well posedness since $t_s=10^{-2}$ is rather small but will be large enough to observe norm inflation in this low regularity regime.

 The Student data, displayed in the upper sub-panels and by the brown markers of Figure~\ref{fig:PWP}, reproduce this picture sub-panel by sub-panel: $\SN$ depends very weakly on time, $\Sinf$ saturates at large $N$, and $\Delta \Sinf$ decays following the same power law $N^{s-(\alpha-1/2)}$. This last point is the most significant one, and it is not a coincidence: the rate \eqref{eq:ConvPWPInitCond} at which the truncated data \eqref{eq:WPInit} converge toward \eqref{eq:GoodInit} involves the variance of $\theta$ and nothing else, so that it is insensitive to the decay of its tails. The two laws are therefore expected to converge at the same rate, and they are observed to do so, in both the energy sub-critical and the energy super-critical case. These simulations thus support Conjecture~\ref{conj:HeavyTail1d} and suggest that the sub-Gaussianity condition \eqref{eq:SubGauss} is a requirement of the available proofs rather than of the phenomenon itself.

\subsection{Norm inflation }\label{sec:NormInflationNum}
\begin{figure}[htb]
	\centering
	\includegraphics[width=0.80\linewidth]{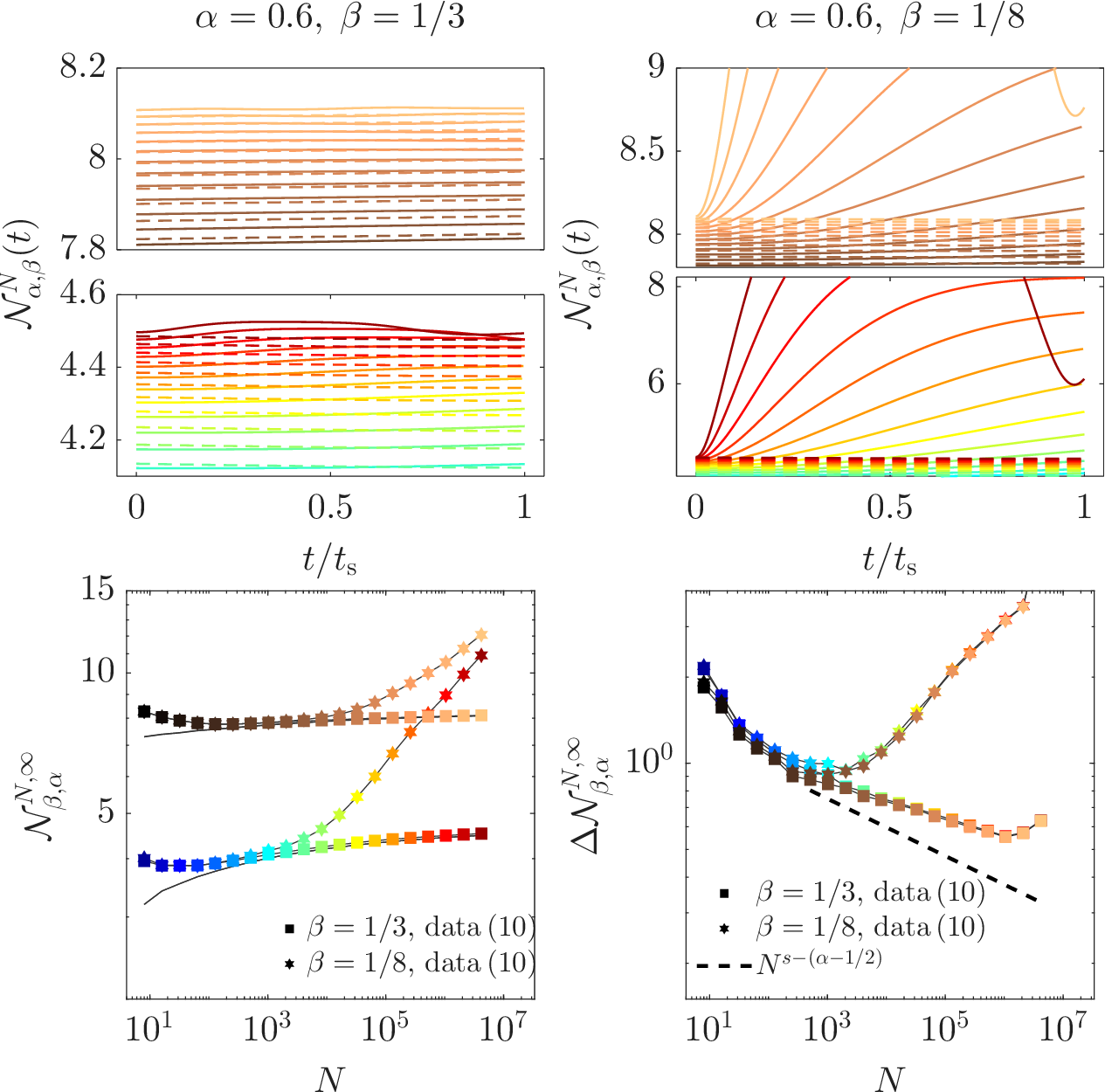}
	\caption{Same figure as for probabilistic well posedness, in the case of norm inflation obtained with the pathological initial data \eqref{eq:PathologicalInit} and $\alpha=0.6$. The upper row shows the time evolution of $\NN$ in the energy sub-critical (left column) and energy super-critical (right column) case in solid colored lines; within each column the upper sub-panel is the Student law $\theta_s$ (black to light brown colormap) and the lower one the Gaussian law $\theta_g$ (blue to red colormap), following the convention of Figure~\ref{fig:PWP}. For reference, we superimpose the results of the upper row of Figure~\ref{fig:PWP} in colored dashed lines, obtained with the initial data \eqref{eq:WPInit}. The left plot of the bottom row shows the evolution of $\Ninf$ with $N$ in both subcritical (squares) and super-critical (stars) cases and for the two laws, the upper group of markers corresponding to the Student law and the lower one to the Gaussian law at small $N$. We also superimpose $\Sinf$ in thin solid lines without any marker to compare to the probabilistic well-posedness case. Finally, the right figure presents the evolution of $\Delta \Ninf$ as a function of $N$ with the same markers as earlier. We also superimpose in thin solid lines the probabilistic well posedness case of Figure \ref{fig:PWP} and, in black dashed line, the power law $N^{s-(\alpha-\frac12)}$ entering the rate of convergence \eqref{eq:ConvPathologicalInitCond} of the initial condition.  }
	\label{fig:NormInflation}
\end{figure}

The results presented in the previous Section \ref{sec:PWPNum} support the extension of probabilistic well-posedness for the one dimensional fractional nonlinear wave equation \eqref{eq:NLWaves} in both energy sub-critical and energy super-critical regimes. The necessity of thinking of well posedness in a probabilistic manner arises from the ill posedness in general of \eqref{eq:NLWaves} for general low regularity data. This ill posedness at low regularity can be seen through norm inflation which we illustrate in this section. Norm inflation is typically proven using a sequence of initial condition concentrating in one point while converging to zero in a given Sobolev space. This is the motivation for the choice of perturbation $p_N^s$ entering in the pathological approximation \eqref{eq:PathologicalInit}. In the regime where the nonlinear behavior of the equation dominates the dispersive one ($\alpha-1/2< 1/2-\beta$), such a concentration of the initial data can lead to instantaneous (in the limit of infinite $N$) growth of Sobolev norms.  In this regime $\alpha-1/2< 1/2-\beta$, we compute the $\Hs{}$ norm of the numerical approximation of the solution of \eqref{eq:NLWaves} with initial condition \eqref{eq:PathologicalInit} that we recall is noted $\NN$ and defined as \eqref{eq:SobolNum}. Similarly to the case of probabilistic well posedness exposed in the previous section, we track the behavior of $\NN$ as a sequence of function and we compare it to $\SN$ obtained with initial data \eqref{eq:WPInit} presented in the previous section. This comparison is displayed in Figure~\ref{fig:NormInflation} where the upper row presents the behavior of $\NN$ (solid colored lines) as a function of time for different values of $N$ and we compare it to $\SN$ (colored dashed lines). As in Figure~\ref{fig:PWP}, the left column is the energy sub-critical case $\beta=1/3$ and the right one the energy super-critical case $\beta=1/8$, and within each column the upper sub-panel is devoted to the Student law $\theta_s$ and the lower one to the Gaussian law $\theta_g$.
In both energy sub-critical ($\beta=1/3$, left column of the upper row) and energy super-critical ($\beta=1/8$, right column of the upper row), and for both laws \eqref{eq:DistribFourier}, the dashed and solid lines are getting closer to each other at time $t=0$, signifying that both approximations \eqref{eq:WPInit} and \eqref{eq:PathologicalInit} are converging toward the same initial condition \eqref{eq:GoodInit}. However, solid and dashed lines are diverging from each other at later times when $N$ increases. This divergence is very striking in the energy super-critical $\beta=1/8$ case but is milder yet still present in the energy sub-critical case $\beta=1/3$. This is indeed expected since in the energy supercritical case $\beta=1/8$, the nonlinearity is stronger (compared to dispersion) than in the subcritical case. In addition, we note that solid lines $\NN$ are diverging sooner and sooner from dashed lines $\SN$ as $N$ gets larger, a signature of the norm inflation phenomenon, for which we expect that 
$$ \lim_{N\to \infty} \sup_{t \in [0,T]} \NN=+\infty,  $$  
for any $T>0$. In other word, the $\Hs{}$ norm of the solution becomes instantaneously infinite in the infinite $N$ limit while still converging at initial time. This divergence in $C([0,T],\Hs{})$ is quantified in the left panel of the bottom row of Figure~\ref{fig:NormInflation} where we show $\Ninf$ as a function of $N$ with colored markers whose color corresponds to the value of $N\in\lbrace 2^3, \dots , 2^{23}\rbrace$ for both energy sub-critical (squares) and energy super-critical (stars) and for the two laws \eqref{eq:DistribFourier}. In addition of these markers, we add in thin solid lines the results presented in Figure~\ref{fig:PWP} for comparison. The results from the energy super-critical $\beta=1/8$ case show a clear divergence from the simulations with initial condition \eqref{eq:WPInit} while this divergence is milder in the energy sub-critical case which is coherent with the heuristic argument that in the energy super-critical case, the nonlinearity dominates the dispersion. Quantifying the growth rate of $\Ninf$ with $N$ as a function of $\alpha$ and $\beta$ would be interesting in itself. However, it would require more resolved numerical simulations to be able to reach higher values of $N$ and therefore being able to fit the behavior at large $N$, which we cannot do with our simulations. Finally, the right panel of the bottom row of Figure~\ref{fig:NormInflation} shows the behavior of $\Delta\Ninf$ as a function of $N$ (colored markers) which is compared to the results obtained in the probabilistic well-posedness regime in thin solid lines. Complementary to the unboundedness in $C([0,T],\Hs{})$, the bottom right panel shows that in the energy super-critical case $\beta=1/8$ the quantity $\Delta \Ninf$ stops decreasing and starts growing when $N$ gets large, indeed signaling that $(w_N)_{N\in \mathbb{N}}$ does not converge  in $C([0,T],\Hs{})$ as expected from its unboundedness. In the energy sub-critical case $\beta=1/3$, $\Delta \Ninf$ still decreases with $N$ but at a rate clearly slower than the power law $N^{s-(\alpha-1/2)}$ observed in Figure~\ref{fig:PWP}, consistently with the merely logarithmic convergence \eqref{eq:ConvPathologicalInitCond} of the pathological approximation and with the milder norm inflation observed in this regime.

Heavy tails do not affect this mechanism. The Student and the Gaussian data of Figure~\ref{fig:NormInflation} follow one another closely: $\Ninf$ saturates for $\beta=1/3$ and grows without bound for $\beta=1/8$, and $\Delta \Ninf$ fails to converge to zero for $\beta=1/8$, in both cases. This is consistent with the fact that norm inflation is driven by the deterministic perturbation $p_N^s$, whose size in $H^s(\T)$ does not depend on the law of $g_\omega$ and $h_\omega$; heavy tails of the Fourier coefficients therefore neither cure nor aggravate the ill posedness caused by the pathological approximation \eqref{eq:PathologicalInit}.

\subsection{The deterministic well posedness regime: a numerical fail-safe }\label{sec:LWPNum}

Given the singular nature of the perturbation $p_N^s$, one might see the divergences exposed in the previous section as a purely numerical effect obtained from a lack of robustness of the numerical scheme. In order to rule out this possibility, we also explore the regime of parameters for which we expect deterministic well-posedness in $\Hs{}$  to hold, namely $1/2-\beta <s< \alpha -1/2$. It means that the solution at positive times depends continuously on the initial data and converges in $C([0,t_s],\Hs{})$. In particular, both sequences $(u_N)_{N\in \mathbb{N}}$ and $(w_N)_{N\in \mathbb{N}}$ obtained from the initial approximations \eqref{eq:WPInit} and \eqref{eq:PathologicalInit} should be bounded and converge in $C([0,t_s],\Hs{})$ to the same solution of \eqref{eq:NLWaves} as $N$ gets large. Since this fail-safe is meant to validate the numerical procedure as a whole, we perform it for the two laws \eqref{eq:DistribFourier}, so that the simulations reported in Sections~\ref{sec:PWPNum} and \ref{sec:NormInflationNum} for the heavy tailed law are subjected to the same test as the Gaussian ones. 
\begin{figure}[htb]
	\centering
	\includegraphics[width=0.80\linewidth]{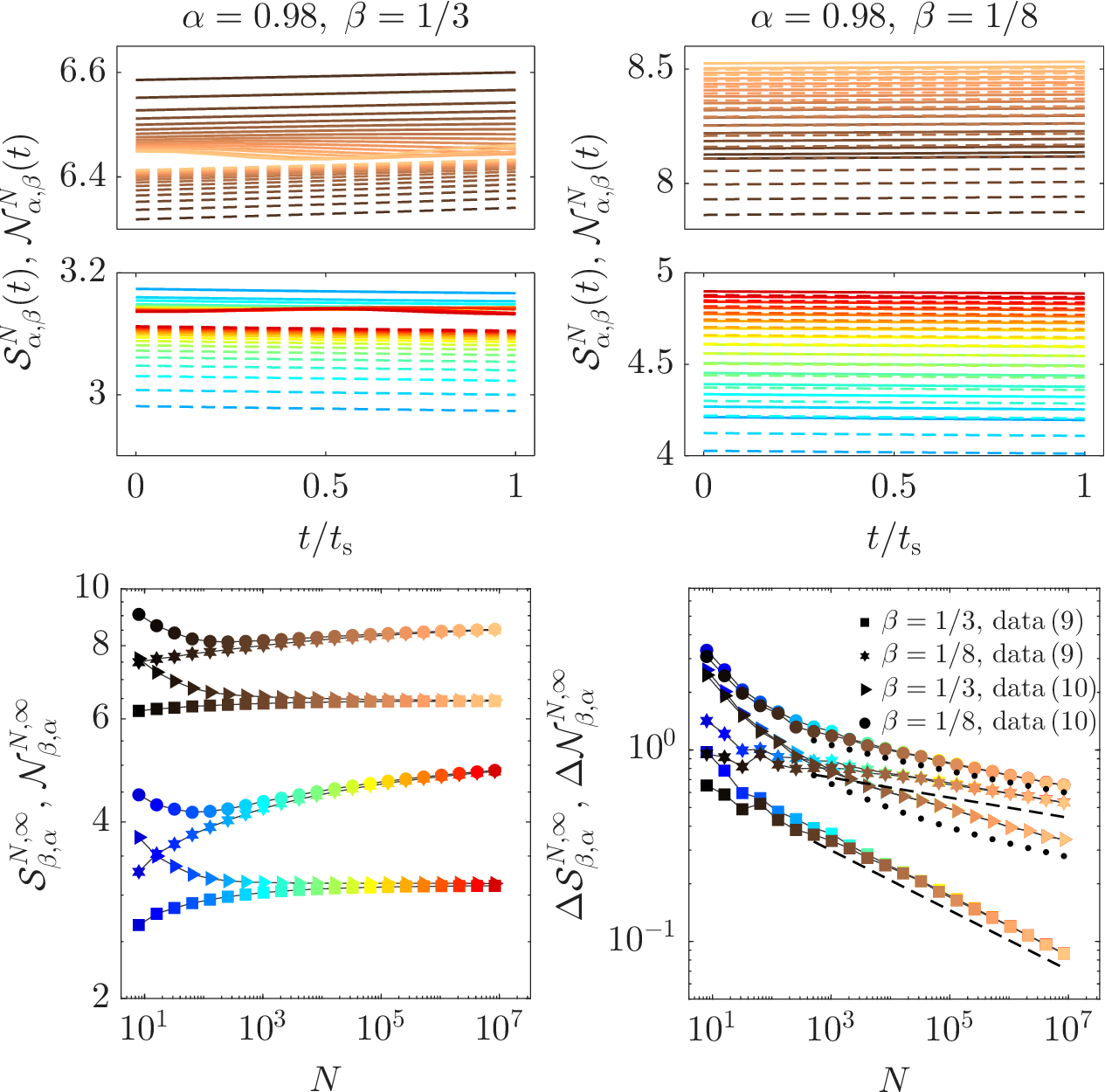}
	\caption{Evolution of the norms in the local well-posedness regime $1/2-\beta<\alpha-1/2$ (here $\alpha=0.98$) for which we expect the solutions obtained from the approximations \eqref{eq:WPInit} and \eqref{eq:PathologicalInit} to converge toward the same unique solution of \eqref{eq:NLWaves}. The upper row shows the behavior of $\NN$ in solid colored lines and $\SN$ in dashed colored lines in the energy subcritical (left column) and supercritical (right column) cases; within each column the upper sub-panel corresponds to the Student law $\theta_s$ and the lower one to the Gaussian law $\theta_g$, with the color conventions of Figure~\ref{fig:PWP}. On the bottom row, the left figure shows the evolution with $N$ of $\Ninf$ and $\Sinf$ for both values of $\beta$ and for the two laws, showing that $\Ninf$ and $\Sinf$ both converge toward the same finite value; the upper group of markers corresponds to the Student law and the lower one to the Gaussian law. The right plot of the bottom row shows the $N$ dependence of $\Delta \Ninf$ and $\Delta \Sinf$, showing convergence to zero at the same rate as the initial condition for both values of $\beta$ and both laws. We also add black dashed and dotted lines corresponding to the rates of convergence of the approximation of the initial data \eqref{eq:ConvPWPInitCond} (dashed lines) and \eqref{eq:ConvPathologicalInitCond} (dotted lines). The slower than power law decay observed can thus be attributed to the logarithmic convergence of \eqref{eq:PathologicalInit} toward \eqref{eq:GoodInit}.  }\label{fig:LWP}
\end{figure}
The results of numerical simulations shown in Figure~\ref{fig:LWP} indeed advocate for local well-posedness for high enough Sobolev regularity ($1/2-\beta <s< \alpha-1/2$) of the initial data, in both sub-critical and super-critical cases. In the upper row of Figure~\ref{fig:LWP}, we display the time dependence of the Sobolev norms of the solutions obtained from approximation \eqref{eq:WPInit} (dashed lines) and \eqref{eq:PathologicalInit} (solid lines). The four sub-panels of the upper row indeed show that dashed colored lines and solid colored lines corresponding respectively to $\SN$ and $\NN$ converge toward the same bounded function of time as $N$ increases, for both laws \eqref{eq:DistribFourier}. This behavior is supported by the left panel of the bottom row of Figure~\ref{fig:LWP} comparing the behaviors of $\Sinf$ and $\Ninf$ for both $\beta=1/8$ and $\beta=1/3$ showing that indeed $\Sinf$ and $\Ninf$ converge with $N$ toward the same finite value. This is in stark contrast with the results of the previous section, showing the divergence of $ \Ninf$ in the large $N$ limit. Finally, the right panel of the bottom row confirms the convergence in $C([0,t_s],\Hs{})$ through the lens of $\Delta \Sinf$ and $\Delta \Ninf$ exhibiting a convergence to zero at the same rate as the one expected for the approximations of the initial condition (black dashed lines). Indeed, the slowest convergence rate displayed in this panel is the one for $\beta=1/3$ and initial approximation \eqref{eq:PathologicalInit} which shows only a logarithmic convergence to zero as shown by the black dotted line, corresponding to \eqref{eq:ConvPathologicalInitCond}. This logarithmic convergence to zero should also be observed in the case of $(\alpha,\beta)= (0.98,1/8)$ displayed by colored circles in Figure~\ref{fig:LWP} but for the values of $N$ we can explore, the power law behavior still dominates on the logarithmic one. As it can be noted in the right panel of the bottom of Figure~\ref{fig:LWP}, the convergence rate to zero depends on the value of $\beta$ while it was not the case in previous Sections~\ref{sec:PWPNum}, \ref{sec:NormInflationNum}. This is caused by our choice $\gamma = 1/2(1 + (1/2-\beta)/(\alpha-1/2))$ which depends on $\beta$ while we chose $\gamma=1/3$ independent of $\beta$ in Sections~\ref{sec:PWPNum} and \ref{sec:NormInflationNum}. The fail-safe is passed by the heavy tailed simulations as well: for each value of $\beta$ and each of the two laws, $\Sinf$ and $\Ninf$ converge toward a common finite value and $\Delta \Sinf$, $\Delta \Ninf$ converge to zero at the predicted rates. This rules out the possibility that the agreement between the two laws observed in Sections~\ref{sec:PWPNum} and \ref{sec:NormInflationNum} is a numerical artifact caused by the occasional large Fourier coefficients produced by $\theta_s$.

\subsection{Numerical conservation of the Hamiltonian}\label{sec:HamilCons}

\begin{figure}
	\centering
	\includegraphics[width=0.45\linewidth]{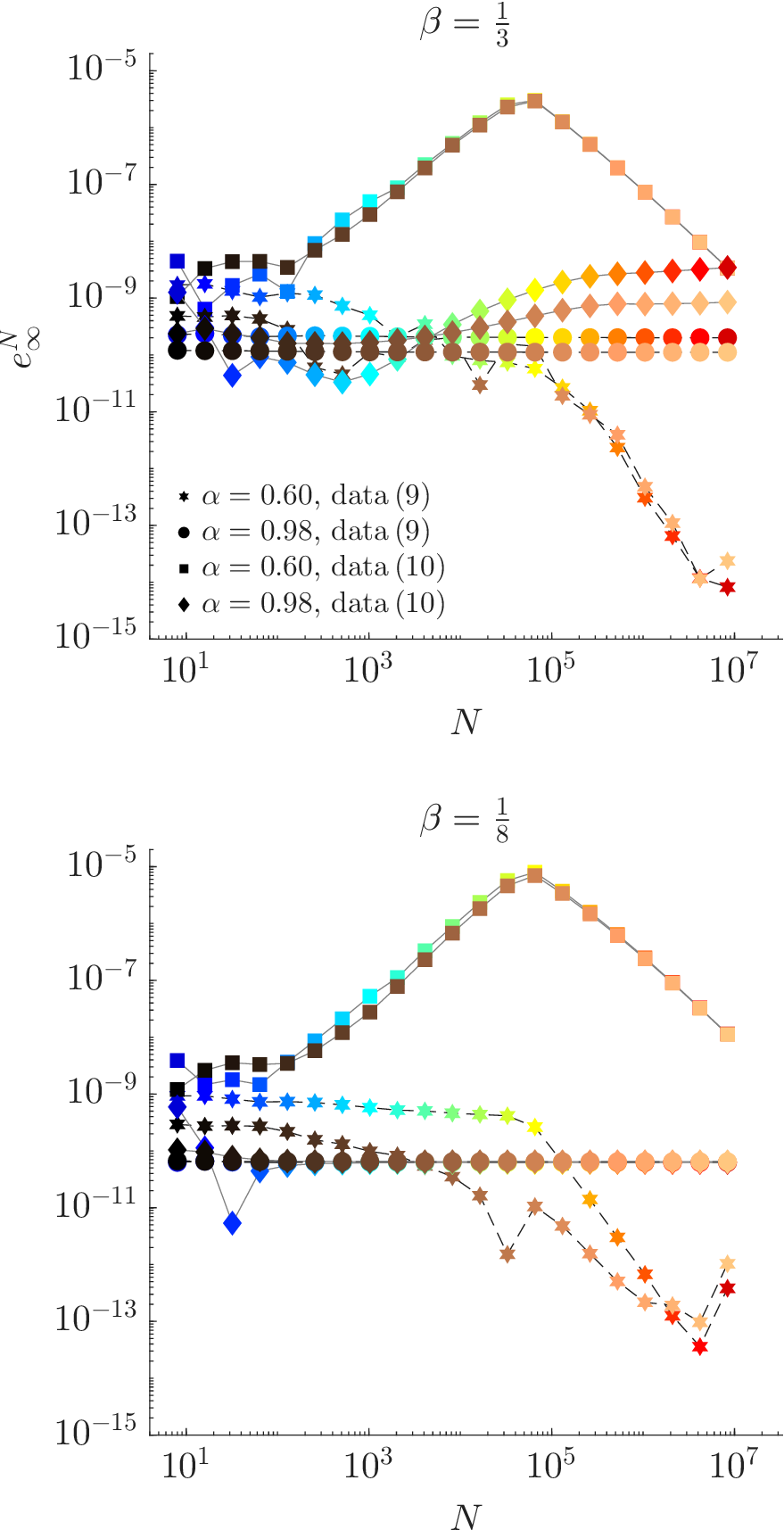} 
	\includegraphics[width=0.45\linewidth]{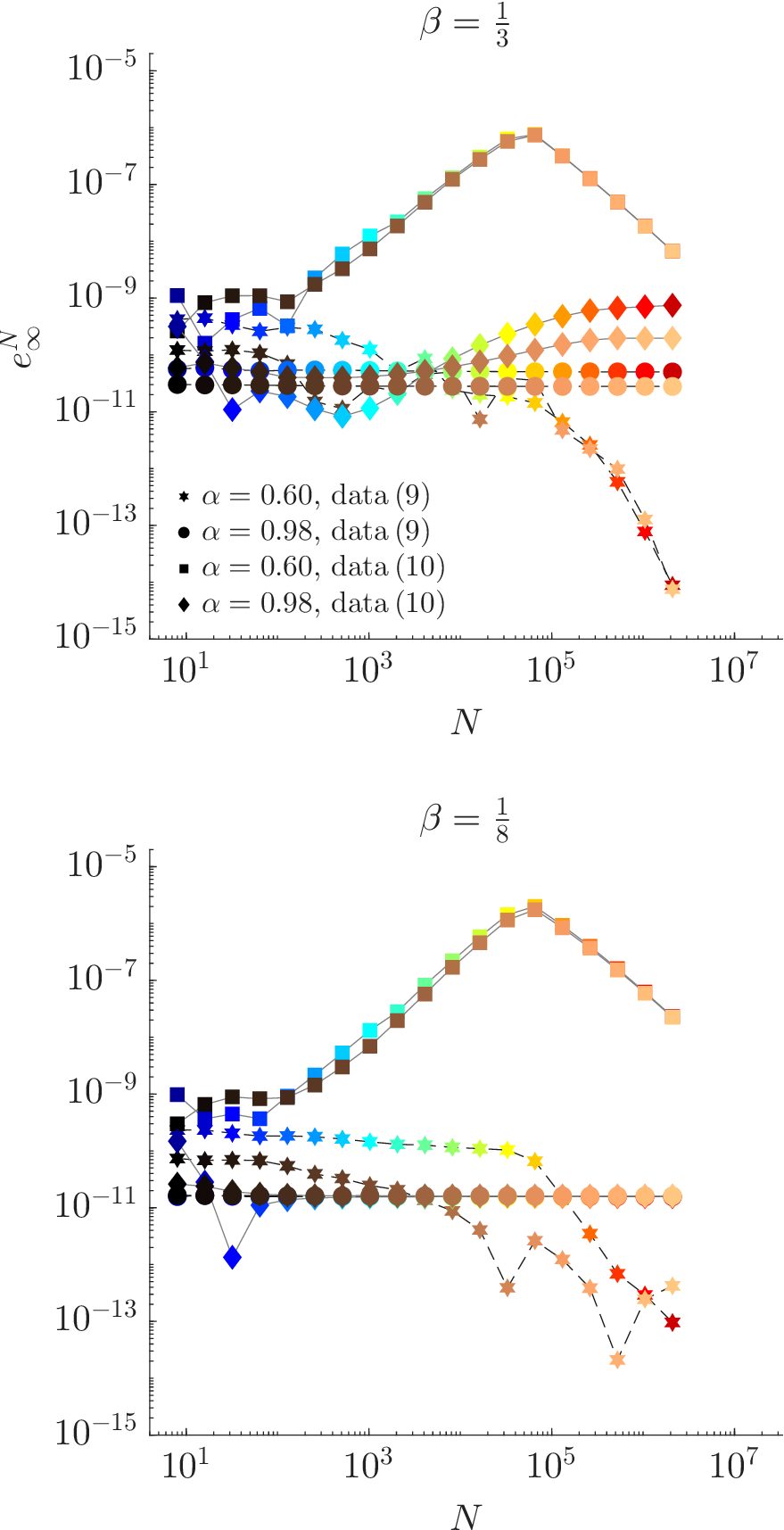} 
	\caption{ Relative error \eqref{eq:NumError} on the conservation of the discretized Hamiltonian. The left column corresponds to the discretization parameters $\tau_N$, $M$ used throughout the paper, and the right column to the refined simulations with $\tau_N/2$, $2M$. In each column, the upper panel shows the error for $\beta=1/3$ and the lower one for $\beta=1/8$. Solid lines correspond to the pathological approximation \eqref{eq:PathologicalInit} while dashed lines are for the unperturbed case \eqref{eq:WPInit}. In each case, circle and diamond markers are used for $\alpha=0.98$ while square and star markers are for $\alpha=0.6$. Both laws \eqref{eq:DistribFourier} are displayed, with the color conventions of Figure~\ref{fig:PWP}: blue to red markers for the Gaussian law $\theta_g$ and black to light brown markers for the Student law $\theta_s$. Comparing the two columns shows that refining the discretization lowers the error everywhere without changing its qualitative dependence on $N$, $\alpha$ and $\beta$.   }
	\label{fig:ErrorHamil}
\end{figure}

In order to probe fine properties of the nonlinear wave equation \eqref{eq:NLWaves}, we have to make sure that the numerical methods are reliable. The only observable we can control here is given by the Hamiltonian structure of \eqref{eq:NLWaves} and should be preserved by the numerical integration in the limit of infinite spatial and temporal resolution. Indeed, for any finite $N$ the solution of \eqref{eq:NLWaves} is smooth, ensuring the convergence of the numerical scheme and also the finiteness of the Hamiltonian. However, we need to make sure that the control of the error remains satisfying in the limit of large $N$. As we explained earlier, we chose empirically a $N$ dependent time step $\tau_N \simeq (N^{1/2-s}/\log N)^ {-3} $. In Figure~\ref{fig:ErrorHamil}, we present the relative error \eqref{eq:NumError} on the Hamiltonian conservation accumulated through time integration. In this figure, each dot corresponds to one of the numerical simulations presented in the previous sections, covering the cases $\alpha=0.6,\, 0.98$ and $\beta=1/3, \, 1/8$, both approximations \eqref{eq:WPInit} and \eqref{eq:PathologicalInit} and both laws \eqref{eq:DistribFourier} of the Fourier coefficients. The right column of Figure~\ref{fig:ErrorHamil} presents the same results for a refined time step $\tau_N/2$ and a refined spatial resolution $2M$, with $N$ ranging up to $2^{21}$, so that the effect of the discretization can be read off by comparing the two columns panel by panel.

For each numerical experiment the error \eqref{eq:NumError} on the conservation of the Hamiltonian is well controlled as shown in Figure~\ref{fig:ErrorHamil}. The relative error remains below $10^{-5}$ in every configuration we explored, which is several orders of magnitude smaller than the variations of the Sobolev norms reported in the previous sections. It is by far the largest for the low regularity, perturbed configuration $\alpha=0.6$ with the approximation \eqref{eq:PathologicalInit}, for which it grows with $N$ up to $N \simeq 5 \times 10^4$ and decreases beyond. This non-monotonicity reflects the crossover in \eqref{eq:TimeStep} between the dispersive constraint $\tau_d/5$, which is the binding one at small $N$ while the solution becomes increasingly violent, and the nonlinear constraint, which takes over at large $N$ and forces $\tau_N$ to vanish rapidly. For the regular data $\alpha=0.98$, and for the unperturbed approximation \eqref{eq:WPInit}, the error is uniformly bounded by $10^{-9}$ and shows no tendency to grow with $N$. Let us stress that the markers associated with the Student law show no systematic departure from the Gaussian ones, the difference between the two laws being small compared with the spread across $\alpha$ and across the two approximations. In particular, the occasional large Fourier coefficients produced by the heavy tailed law $\theta_s$ do not degrade the conservation of the Hamiltonian, which confirms that the time step \eqref{eq:TimeStep}, calibrated using the variance of $\theta$ only, is also adequate outside of the sub-Gaussian framework. 

Comparing the two columns of Figure~\ref{fig:ErrorHamil} shows that the refinement $\tau_N/2$, $2M$ lowers the error in every panel while leaving its shape unchanged. In particular the maximum of the error, reached at $N \simeq 5\times 10^4$ for the perturbed low regularity configuration, is divided by about four, which is consistent with the second order accuracy in time of the scheme \eqref{eq:NumScheme}; the plateaus attained by the regular data $\alpha=0.98$ are lowered by a comparable factor. For reasons of numerical resources the refined simulations only reach $N=2^{21}$, whereas the reference ones go up to $N=2^{23}$. All in all, the error on the Hamiltonian is well controlled over all the explored regimes. In the infinite $N$ limit, the solution looses its smoothness and the Hamiltonian \eqref{eq:HamiltonianNLWave} need not remain finite. Monitoring the conservation of a diverging quantity might seem unfit for the control of numerical error. However, we recall that both the time step and the spatial step vanish in the large $N$ limit and that convergence of the numerical scheme should also break down. The numerical approach to probabilistic well-posedness and norm inflation must therefore break down in the limit irrespectively of the numerical method. For these reasons, the numerical simulations must be seen as a finite $N$ illustration of the limiting behavior.
\section{Conclusion and perspectives}\label{sec:conclu}
The cubic defocusing nonlinear wave equation is generally ill posed in Sobolev spaces with low regularity $s<s_c$. This ill posedness can manifest through norm inflation, an instantaneous growth of the Sobolev norm of the initial data. However, well-posedness can be recovered on a probabilistic level using random initial data whose Fourier coefficients are sub-Gaussian, together with a given--Fourier truncation--approximation of such initial data. Here, we studied numerically the one dimensional and fractional cubic defocusing wave equation allowing us to explore energy sub-critical and super-critical regimes. Our numerical results show that it is possible to illustrate numerically very fine properties of the nonlinear wave equation such as probabilistic well posedness. We also made sure that for high enough regularity of the initial data, the solution converges independently of the approximation of the initial condition in $\Hs{}$ as expected from local well posedness.

The main outcome of this work concerns the role of the sub-Gaussianity condition \eqref{eq:SubGauss}. As recalled in Remark~\ref{rk:RoleSubGauss}, this condition is used in an essential way in the proof of Theorem~\ref{thm:ProbaWP}, through large deviation estimates whose Gaussian tails are what makes the almost sure argument work; laws with power law tails are, in this respect, out of reach of the currently available methods. We therefore replaced the Gaussian law by a Student law with three degrees of freedom, which is centered with unit variance but violates \eqref{eq:SubGauss}, and repeated the whole numerical study. All the phenomena described above are reproduced without qualitative change: boundedness and convergence in $C([0,t_s],\Hs{})$ of the Fourier truncated sequence, norm inflation for the pathological approximation \eqref{eq:PathologicalInit}, and convergence of both approximations toward the same solution in the deterministic well posedness regime. Moreover the measured convergence rate is unchanged, which is consistent with the fact that the rate \eqref{eq:ConvPWPInitCond} at which the truncated initial data converge involves the variance of the law and nothing else.

These results are the numerical evidence we can offer in favour of Conjecture~\ref{conj:HeavyTail1d}, and indirectly of Conjecture~\ref{conj:HeavyTail}, according to which the conclusion of Theorem~\ref{thm:ProbaWP} should persist when \eqref{eq:SubGauss} is replaced by the sole requirement that $\theta$ be centered with finite variance. To the best of our knowledge this is a completely open problem, and we hope that the present study will motivate work in that direction. Two questions seem particularly natural to us. The first one is how many moments of $\theta$ are actually needed, if any beyond the second; our simulations use a law with moments of order $p$ for every $p<3$ and say nothing about heavier tails. The second one is whether the global in time conclusion of Theorem~\ref{thm:ProbaWP} survives outside of the sub-Gaussian framework, where the energy of the random initial data is itself heavy tailed; our simulations probe the short time $t_s=10^{-2}$ and bear on the local in time statement only.

Let us close with a perspective of a different nature. The failure of a sequence of regularized problems to select a single limit is not specific to the Fourier truncation considered here. It is also at the heart of spontaneous stochasticity, where a vanishing viscosity, noise or regularization scale fails to restore uniqueness and the singular limit is described by a probability law on solutions rather than by a single one, see \cite{ruffenach2026measure} and the references therein. This phenomenon has recently been reported in a dispersive setting, for a focusing Majda--McLaughlin--Tabak equation undergoing finite time wave collapse, where vanishing perturbations of the regularization parameter or of the initial condition survive the singular limit and generate a finite post-blowup uncertainty \cite{ruffenach2026collapse}. This is of direct relevance to the regimes probed in the present paper. On the one hand, the energy super-critical case $\beta<\beta_c$, in which finite time blowup of smooth solutions is not excluded, is a natural candidate for such a behavior. On the other hand, the pathological approximation \eqref{eq:PathologicalInit} provides a sequence which is precisely observed not to converge in Section~\ref{sec:NormInflationNum}, and it would be interesting to know whether the associated sequence $(w_N)_{N \in \mathbb{N}}$ admits a description in terms of a limiting law rather than of a limiting solution.

\section*{Acknowledgments}
The authors thank ENS de Lyon for funding and for numerical infrastructures. The authors also thank Erwan Faou for his helpful comments and suggestions concerning the elaboration of this manuscript.

\bibliographystyle{unsrt}

\end{document}